\documentclass[preprint,sort&compress]{elsarticle}
\usepackage[margin=1in]{geometry}
\usepackage{graphicx}
\usepackage{pstool}
\usepackage{wrapfig}
\usepackage{caption}
\usepackage{subcaption}
\usepackage[section]{placeins} % Defines \FloatBarrier which does not allow
\usepackage{fancyhdr} % Required for custom headers
\usepackage{lastpage} % Required to determine the last page for the footer
\usepackage{extramarks} % Required for headers and footers

\usepackage{xcolor} % Colors
\usepackage{enumerate} % Additional options for \enumerate{}
\usepackage{paralist} % Additional options for lists; inline lists
\usepackage{amsmath, amsthm, amssymb, mathtools}
\usepackage{mathabx, pifont, stmaryrd} % More math symbols, fonts
\usepackage[explicit]{titlesec} % Control over section/subsection styles
\usepackage{etoolbox} % Boolean variables and more
\usepackage{bibentry} % Bibliography entries anywhere in doc
\makeatletter\let\saved@bibitem\@bibitem\makeatother % Do not remove if bibentry
\usepackage[colorlinks, bookmarksopen, bookmarksnumbered,
            citecolor=red,urlcolor=red]{hyperref} % Hyperlinks
\makeatletter\let\@bibitem\saved@bibitem\makeatother % Do not remove if bibentry
\usepackage{algorithm}
\usepackage{algorithmic}

\usepackage{bm} % required for: \bm{}
\usepackage{amsfonts} % required for: \mathbb{}, \mathcal{}, ...

\newcommand{\func}[3]{\ensuremath{#1 : #2 \rightarrow #3}}

\newcommand{\norm}[1]{\ensuremath{\left\| #1 \right\|}}

\newcommand{\optunc}[2]{\underset{#1}{\text{minimize}} ~~ #2}

\newcommand{\pder}[2]{\ensuremath{\frac{\partial #1}{\partial #2}}}

\newcommand{\oder}[2]{\ensuremath{\frac{d#1}{d#2}}}
\newcommand{\Dcal}{\ensuremath{\mathcal{D}}}

\newcommand{\Ncal}{\ensuremath{\mathcal{N}}}
\newcommand{\Ocal}{\ensuremath{\mathcal{O}}}

\newcommand{\Tcal}{\ensuremath{\mathcal{T}}}

\newcommand{\Vcal}{\ensuremath{\mathcal{V}}}
\newcommand{\Wcal}{\ensuremath{\mathcal{W}}}

\newcommand{\Rbb}{\ensuremath{\mathbb{R} }}

\newcommand\Abm{{\ensuremath{\bm{A}}}}

\newcommand\Ibm{{\ensuremath{\bm{I}}}}

\newcommand\Mbm{{\ensuremath{\bm{M}}}}

\newcommand\Pbm{{\ensuremath{\bm{P}}}}

\newcommand\Xbm{{\ensuremath{\bm{X}}}}

\newcommand\fbm{{\ensuremath{\bm{f}}}}

\newcommand\ubm{{\ensuremath{\bm{u}}}}
\newcommand\vbm{{\ensuremath{\bm{v}}}}
\newcommand\wbm{{\ensuremath{\bm{w}}}}

\newcommand\ybm{{\ensuremath{\bm{y}}}}

\newcommand\mubold{{\ensuremath{\boldsymbol{\mu}}}}

\newcommand\taubold{{\ensuremath{\boldsymbol{\tau}}}}

\newcommand\nubold{{\ensuremath{\boldsymbol{\nu}}}}

\newcommand\xibold{{\ensuremath{\boldsymbol{\xi}}}}

\newcommand\Phibold{{\ensuremath{\boldsymbol{\Phi}}}}
\newcommand\Pibold{{\ensuremath{\boldsymbol{\Pi}}}}

\newcommand\Psibold{{\ensuremath{\boldsymbol{\Psi}}}}
\newcommand\Xibold{{\ensuremath{\boldsymbol{\Xi}}}}

\newcommand\onebold{\ensuremath{\mathbf{1}}}

\usepackage{tikz}
\usepackage{pgfplots}
\usepackage{pgfplotstable, filecontents, booktabs}
\pgfplotsset{compat=1.9}

\usetikzlibrary{pgfplots.groupplots}
\usepgfplotslibrary{fillbetween}
\usetikzlibrary{calc,fit,matrix,arrows,automata,positioning,shapes}
\usetikzlibrary{arrows.meta}

\pgfplotsset{select coords between index/.style 2 args={
    x filter/.code={
        \ifnum\coordindex<#1\fi
        \ifnum\coordindex>#2\fi
    }
}}

\tikzset{
 invisible/.style={opacity=0},
 visible on/.style={alt={#1{}{invisible}}},
 alt/.code args={<#1>#2#3}{%
   \alt<#1>{\pgfkeysalso{#2}}{\pgfkeysalso{#3}}
 },
}

\newbool{fastcompile}
\setbool{fastcompile}{false}

\usepackage{lineno}
\usepackage{todonotes}

\journal{Elsevier}

\begin{document}

\begin{frontmatter}
%\tableofcontents
% \listoftodos

%% Title, authors and addresses
\title{Non-intrusive model reduction of large-scale, nonlinear dynamical systems using deep learning}

%% use the tnoteref command within \title for footnotes;
%% use the tnotetext command for the associated footnote;
%% use the fnref command within \author or \address for footnotes;
%% use the fntext command for the associated footnote;
%% use the corref command within \author for corresponding author footnotes;
%% use the cortext command for the associated footnote;
%% use the ead command for the email address,
%% and the form \ead[url] for the home page:
%%
%% \title{Title\tnoteref{label1}}
%% \tnotetext[label1]{}
%% \author{Name\corref{cor1}\fnref{label2}}
%% \ead{email address}
%% \ead[url]{home page}
%% \fntext[label2]{}
%% \cortext[cor1]{}
%% \address{Address\fnref{label3}}
%% \fntext[label3]{}

%% use optional labels to link authors explicitly to addresses:
%% \author[label1,label2]{<author name>}
%% \address[label1]{<address>}
%% \address[label2]{<address>}

\author[ndAME,ndCICS]{Han Gao}
\author[ndAME,ndCICS]{Jian-Xun Wang}
\author[ndAME,ndCICS]{Matthew J. Zahr\corref{corxh}}

\address[ndAME]{Department of Aerospace and Mechanical Engineering, University of Notre Dame, Notre Dame, IN}
\address[ndCICS]{Center for Informatics and Computational Science, University of Notre Dame, Notre Dame, IN}

\cortext[corxh]{Corresponding author. Tel: +1 574-631-1298}
\ead{mzahr@nd.edu}

\begin{abstract}
Projection-based model reduction has become a popular approach to
reduce the cost associated with integrating large-scale dynamical
systems so they can be used in many-query settings such as
optimization and uncertainty quantification. For nonlinear systems,
significant cost reduction is only possible with an additional layer of
approximation to reduce the computational bottleneck of evaluating
the projected nonlinear terms. Prevailing methods to approximate
the nonlinear terms are code intrusive, potentially requiring
years of development time to integrate into an existing codebase,
and have been known to lack parametric robustness.

This work develops a non-intrusive method to efficiently and
accurately approximate the expensive nonlinear terms that arise 
in reduced nonlinear dynamical system using deep neural networks. 
The neural network is trained using only the simulation data 
used to construct the reduced basis and evaluations of the 
nonlinear terms at these snapshots. Once trained, the neural network-based 
reduced-order model only requires forward and backward propagation through 
the network to evaluate the nonlinear term and its derivative, which are used
to integrate the reduced dynamical system at a new parameter
configuration. We provide two numerical experiments---the dynamical
systems result from the semi-discretization of parametrized, nonlinear,
hyperbolic partial differential equations---that show, in addition to
non-intrusivity, the proposed approach provides more stable and
accurate approximations to each dynamical system across a large number
of training and testing points than the popular empirical interpolation
method.
\end{abstract}

\begin{keyword}
%% keywords here, in the form: keyword \sep keyword
  nonlinear model reduction \sep non-intrusive hyperreduction \sep deep learning
\end{keyword}
\end{frontmatter}

%\linenumbers
% \input{highlights}

\section{Introduction}
\label{sec:intro}
Numerical simulations have made an undeniable impact on the fields of science,
engineering, and medicine due to the possibility to study or analyze a
physical system in a virtual setting. However, modeling and simulation of most
practical systems is a computationally expensive endeavor, usually requiring
days on a supercomputer, essentially limiting users to a few runs. However,
the true power of computational physics lies in many-query analyses, e.g.,
optimization and uncertainty quantification (UQ), which require simulations
at a large number of parameter configurations. For example, optimization
problems are ubiquitous in science and engineering and their solutions can
lead to systems with unprecedented efficiency (e.g., energetically optimal
flapping flight \cite{zahr2016dgopt, zahr2016periodic,
                      wang2017energetically, zahr2018flapopt}),
help gain insight to physical phenomena, or determine properties of a system
from sparse, noisy observations of the solution (inverse problems). To enable
these types of many-query analyses on important problems, the issue of high
computational cost of a single simulation must be addressed. Reduced-order
models (ROMs) offer a promising means to do so.

The number of degrees of freedom (DOF) of a dynamical system is dramatically
reduced in ROMs by constraining the dynamics to evolve in a very
low-dimensional subspace, spanned by a set of reduced basis (RB) functions.
These basis functions are usually learned through training data, i.e.,
solution snapshots obtained from high-dimensional model (HDM) simulations.
Although the number of DOF can be significantly reduced by RB projection, for
nonlinear problems, the speedup of the standard ROMs is often marginal due
to a large number of high-dimensional operations involved in the evaluation
of nonlinear terms in the ROM system. Therefore, an additional step, usually
called \emph{hyperreduction}, must be taken to efficiently approximate
nonlinear terms. Most existing hyperreduction techniques, e.g., empirical
interpolation method (EIM)~\cite{barrault2004empirical} and its discrete
variant discrete EIM (DEIM)~\cite{chaturantabut2010nonlinear}, approximate
the HDM nonlinear velocity function using a low-dimensional subspace as well,
which provides a satisfactory tool to deal with nonlinear PDE systems in an
efficient way. Massive speedups can be gained by hyperreduced ROMs in many
cases where the solution of a dynamical system and its nonlinear terms are
well-approximated in low-dimensional subspaces, including non-parametric
problems (reproduce training data), linear elliptic PDEs, or problems with limited
parameter variations
\cite{carlberg2011gnat, farhat2015ecsw, zahr2017multilevel}.

Despite the tremendous promise of ROMs, standard hyperreduction techniques
often struggle to provide a stable and accurate approximation of nonlinear
terms and present notable limitations in parametric settings~\cite{washabaugh2016phd}.
This lack of parametric robustness remains the main roadblock of ROMs being
the technology that enables large-scale many-query analyses which inherently
involve parametric variations. In addition, standard hyperreduction techniques
are code intrusive \cite{chaturantabut2010nonlinear}, usually requiring
hundreds of person-hours to implement, which poses great challenges to
leveraging existing open-source or commercial legacy codes for computational
mechanics. Therefore, there is an increasing interest in developing
non-intrusive or weakly-intrusive ROMs without the need for access to HDM
operators to achieve better robustness and
stability~\cite{audouze2013nonintrusive,peherstorfer2016data,chen2018greedy,yeh2018common}. 
For example, Audouze et al.~\cite{audouze2013nonintrusive}
proposed a non-intrusive proper orthogonal decomposition (POD)-based ROM,
where the Galerkin projection is bypassed by using a radial basis
regression to estimate the RB coefficients directly and does not
require hyperreduction for efficiency. Peherstofer and
Wilcox~\cite{peherstorfer2016data} proposed a data-driven operator inference
approach based on least square fitting to establish a non-intrusive
projection-based ROM framework. Reduced-order models based on piecewise
polynomial approximation of the dynamical system in state space
\cite{rewienski2003trajectory, dong2003piecewise} are popular in subsurface
flow \cite{cardoso2010linearized, he2011enhanced, trehan2016trajectory} and
electrical engineering applications
\cite{rewienski2003trajectory, vasilyev2006macromodel}, but are difficult to
train since they sensitive to the choice of expansion points
\cite{zahr2010mortestbed}.
%\todo{Maybe review some non-intrusive ROM here, like TPWL?}

Recent advances in scientific machine learning (SciML) has been receiving
increased attention in the computational modeling
community~\cite{lee2018basic,brunton2019machine,wang2017physics,
                shanahan2018machine,brunton2019data} and offers new
opportunities to develop more efficient and accurate reduced-order models.
A growing amount of research using machine learning, particularly deep
learning techniques, for model reduction has been witnessed most recently.
Specifically, a majority of these studies focused on learning the closure
model (or discrepancy terms) of projection-based ROMs from data to improve
their predictive
accuracy~\cite{trehan2017error,san2018machine,san2018neural,pan2018data,
               wan2018data,mohebujjaman2019physically,maulik2019time,
               mou2019data}. San and Maulik~\cite{san2018machine,san2018neural}
employed feedforward neural networks (NNs) to build the ROM closures, with
which the performance can be notably improved. Pan and
Duraisamy~\cite{pan2018data} modeled the truncated dynamics in a data-driven
way using sparse regression and neural networks. In a similar vein,
Mohebujjaman et al.~\cite{mohebujjaman2019physically} proposed a data-driven
correction ROM (DDC-ROM) framework, which has been tested on a number of fluid
dynamic problems. In addition to building closure models, machine learning has
also been used to construct more representative basis functions for model
reduction. For example, Lee and Carlberg~\cite{lee2018basic} applied deep
convolutional autoencoders to learn the nonlinear manifold defined by the
parametrized dynamical system solution, which was shown to outperform the
linear POD basis. Another important direction of using machine learning in
model reduction is the development of accurate non-intrusive
ROMs~\cite{xiao2015non,kani2017dr,mohan2018deep,hesthaven2018non,wang2018model,
           san2019artificial,regazzoni2019machine,wang2019non,li2019deep,
           lui2019construction,xie2019non,pawar2019deep,jin2019deep,
           santo2019data,maulik2019using}.
In most of these works, the general idea is to use machine learning models to
learn the temporal evolution of the states in the reduced subspace and thus
the Galerkin projection and intrusive hyperreduction are bypassed. For
example, the dynamics in low-dimensional manifold can be learned using a
multi-level perceptrons
(MLP)~\cite{wang2018model,san2019artificial,wang2019non}, deep residual
recurrent neural networks (RNN)~\cite{kani2017dr}, or Long Short Term Memory
(LSTM) based RNN~\cite{mohan2018deep,li2019deep,maulik2019using}.
Most of these approaches are
purely data-driven and equation-free, which makes it difficult to respect the
underlying PDE structure. Moreover, these works are focused on problems
with non-parametric settings~\cite{lui2019construction,maulik2019using}.

In this work, we propose a novel method to approximate the nonlinear terms
arising in projection-based ROMs that uses deep learning (DL) to
overcome the parametric robustness issues and code intrusion of existing
hyperreduction methods. Namely, a deep fully-connected neural network (NN)
will be built to learn the nonlinear velocity function in the ROM equations by
leveraging the same HDM solution data used to construct the POD basis
and the corresponding velocity data.
The deep NN model here is embedded into the standard projection-based ROM
setting and the resulting dynamical system is solved using numerical time
integration. The proposed method is \textit{non-intrusive}
in the sense that the hyperreduced model is a small dynamical system
with velocity function defined by the NN that can be run independently of
the original simulation code once the NN has been trained.
The performance of the proposed method will be compared against
a ROM without hyperreduction and a ROM with the (D)EIM hyperreduction. Note
that this work is focused on parametrized, nonlinear dynamical systems. To the
best of our knowledge, the current work is the first attempt to build a
DL-based hyperreduction for projection-based ROMs in parametric settings. 
%A similar approach was proposed in \cite{maulik2019using} that uses machine
%learning to approximate the reduced velocity function; however, the

The rest of the paper is organized as follows.
Section~\ref{sec:mor} introduces projection-based model reduction for
nonlinear dynamical systems and briefly discusses the popular, intrusive
hyperreduction method (D)EIM. Section~\ref{sec:romnn} presents the
detailed formulation and training procedure for the proposed non-intrusive
reduced-order model that uses deep learning to approximate the reduced
velocity function. Section~\ref{sec:num-exp} compares the accuracy of the
proposed method against classical model reduction and (D)EIM using
two dynamical systems that result from the semi-discretization of nonlinear,
hyperbolic PDEs. Finally, Section~\ref{sec:conclude}  concludes the paper.

%\begin{itemize}
% \item Other methods use FEM \cite{} or DG \cite{} to directly approximate
%       $\fbm_r$, which is a more direct approach. Still code intrusive
%       because rely on partial assembly.
% \item Non-intrusive MOR approaches exist (TPWL/Q) but very limited
%       prediction capabilities and sensitive to choice of linearization
%       points; not considered
% \item NN training and evaluation well-suited for GPU, (D)EIM inherits
%       parallelism from CM code, some of which employ methods ill-suited
%       for GPUs
%\end{itemize}

\section{Classical model reduction of nonlinear dynamical systems}
\label{sec:mor}

\subsection{Large-scale, nonlinear dynamical system}
\label{sec:mor:hdm}
Consider a parametrized, nonlinear dynamical system that we will take to be
the HDM,
\begin{equation} \label{eqn:semidisc}
 \Mbm\oder{\ubm}{t} = \fbm(\ubm,t,\mubold), \qquad
 \ubm(0) = \ubm_0,
\end{equation}
where $\Dcal \subset \Rbb^{N_\mubold}$ is the parameter space, $\func{\ubm}{[0,T] \times \Dcal}{\Rbb^{N_\ubm}}$ is the time- and parameter-dependent state, $\ubm_0 \in \Rbb^{N_\ubm}$ is the initial condition, $\func{\fbm}{\Rbb^{N_\ubm}\times [0, T] \times \Dcal}{\Rbb^{N_\ubm}}$ is the velocity of the nonlinear dynamical system $(\xibold, t, \mubold) \mapsto \fbm(\xibold, t, \mubold)$,
and $\Mbm \in \Rbb^{N_\ubm\times N_\ubm}$ is the mass matrix. To
approximately integrate the system in (\ref{eqn:semidisc}), we
introduce a discretization of the time domain into $N_t$ intervals
with endpoints $\Tcal \coloneqq \{t_0, t_1, \dots, t_{N_t}\}$ such
that $t_0 = 0$, $t_{N_t} = T$, and $t_n < t_{n+1}$ for
$n = 0, \dots, N_t-1$ and use any standard solver, e.g.,
backward differentiation formulas or Runge-Kutta methods, to
approximate the solution at these nodes.
In this work, we assume the system in (\ref{eqn:semidisc}) is
large-scale ($N_\ubm \gg 1$) and computationally intensive to numerically integrate. Of particular interest are dynamical
systems that result from the semi-discretization of partial
differential equations for large, complex systems, e.g., in computational
fluid dynamics it is not uncommon to have semi-discrete models with
$\Ocal(10^8)$ degrees of freedom \cite{carlberg2011gnat, washabaugh2016phd}.

We assume the computational complexity of evaluating the
$N_\ubm$-components of the velocity function
$(\xibold,t,\mubold)\mapsto\fbm(\xibold,t,\mubold)$ is
$\Ocal(g(N_\ubm))$. Furthermore, we assume the complexity
of evaluating the Jacobian matrix
$(\xibold,t,\mubold)\mapsto\pder{\fbm}{\xibold}(\xibold,t,\mubold)$
is also $\Ocal(g(N_\ubm))$. This holds for local discretizations
such as the finite difference, finite element, or finite
volume methods where the Jacobian is sparse with the
number of nonzeros per column and computational complexity
of evaluating each column independent of $N_\ubm$. The
complexity of an entire implicit time step---dominated by
the velocity and Jacobian evaluation and the linear solve
with the Jacobian matrix---is $\Ocal(g(N_\ubm)+N_\ubm^3)$
if a direct solver is used and $\Ocal(g(N_\ubm)+N_\ubm^2)$
if an iterative solver with an effective preconditioner is
used.
%Complexity HDM: $f$ is $\Ocal(g(N_\ubm))$ and $df$ is $\Ocal(g(N_\ubm))$;
%constant consumed in $\Ocal$ notation usually much larger for Jacobian
%evaluation; assume local method such as FD, FV, FE where number of
%non-zero entries in each column of the Jacobian does not scale with
%$N_\ubm$; entire implicit time step $\Ocal(g(N_\ubm)+N_\ubm^3)$ or
%$\Ocal(g(N_\ubm)+N_\ubm^2)$ if direct/iterative solve.
%

\subsection{Projection-based model order reduction}
\label{sec:mor:rom}
The construction of projection-based reduced-order models begins with the
ansatz that the solution of the dynamical system can be well-approximated
in a low-dimensional affine subspace
$\Vcal \coloneqq \{\bar\ubm+\Phibold\ybm \mid \ybm \in \Rbb^{k_\ubm}\}$,
where $\Phibold \in \Rbb^{N_\ubm\times k_\ubm}$ with
$\Phibold^T\Phibold = \Ibm$ denotes the reduced basis of
dimension $k_\ubm \ll N_\ubm$ and $\bar{\ubm} \in \Rbb^{N_\ubm}$
is an affine offset. That is,
\begin{equation} \label{eqn:rom-ansatz}
    \ubm(t, \mubold) \approx \ubm_r(t, \mubold) \coloneqq
    \bar\ubm + \Phibold\ybm(t,\mubold),
\end{equation}
where $\func{\ubm_r}{[0,T] \times \Dcal}{\Vcal}$ is the 
approximation of $\ubm(t,\mubold)$ in the reduced subspace
and $\func{\ybm}{[0,T] \times \Dcal}{\Rbb^{k_\ubm}}$ are the reduced
coordinates of $\ubm_r(t,\mubold)$ corresponding to the basis
$\Phibold$ and offset $\bar\ubm$. The reduced coordinates are defined by
enforcing the subspace approximation (\ref{eqn:rom-ansatz}) in the governing
equation and constraining the resulting system to be orthogonal to a
test subspace $\Wcal$ of dimension $\dim\Wcal = k_\ubm$
\begin{equation} \label{eqn:rom}
    \Mbm_r\oder{\ybm}{t} = \fbm_r(\ybm, t, \mubold), \qquad
    \ybm(0) = \ybm_0,
\end{equation}
where $\Psibold \in \Rbb^{N_\ubm\times k_\ubm}$ with
$\Psibold^T\Psibold = \Ibm$ is a basis for $\Wcal$. The
velocity of the reduced dynamical system is
\begin{equation} \label{eqn:rom-velo}
   \func{\fbm_r}{\Rbb^{k_\ubm}\times [0, T] \times \Dcal}{\Rbb^{k_\ubm}}, \quad
   (\taubold, t, \mubold) \mapsto
   \Psibold^T\fbm(\bar\ubm+\Phibold\taubold, t, \mubold)
\end{equation}
and the reduced mass matrix $\Mbm_r\in\Rbb^{k_\ubm\times k_\ubm}$ and
initial condition for the reduced coordinates
$\ybm_0 \in \Rbb^{k_\ubm}$ are defined as
\begin{equation} \label{eqn:red-mass-ic}
    \Mbm_r \coloneqq \Psibold^T\Mbm\Phibold, \qquad
    \ybm_0 \coloneqq \Phibold^T(\ubm_0-\bar{\ubm}).
\end{equation}
The reduced initial condition $\ybm_0$ is the orthogonal projection of
the initial condition $\ubm_0$ onto the trial subspace.

In this work, we choose the test space to be the same as the
trial space, up to the affine offset, i.e., a Galerkin projection
$\Psibold = \Phibold$. The affine offset is taken to be the initial
condition $\bar\ubm = \ubm_0$ based on the observations in
\cite{amsallem2012nonlinear}. The reduced basis $\Phibold$ is
defined using the method of snapshots and proper orthogonal
decomposition (POD) \cite{sirovich1987turbulence}.
For each parameter in a given training set
$\Xibold_0 \coloneqq \{\mubold_1,\dots,\mubold_{N_s}\} \subset \Dcal$,
we compute the approximate solution on the time discretization $\Tcal$ and agglomerate into a global snapshot matrix
$\Xbm = \begin{bmatrix} \Xbm_1 & \cdots & \Xbm_{N_s}\end{bmatrix} \in \Rbb^{N_\ubm \times N_tN_s}$, where fixed-parameter snapshot matrices
$\Xbm_k \in \Rbb^{N_\ubm\times N_t}$ are defined as
\begin{equation} \label{eqn:snapmat}
    \Xbm_k \coloneqq
    \begin{bmatrix}
     \ubm(t_1, \mubold_k) & \cdots & \ubm(t_{N_t}, \mubold_k)
    \end{bmatrix}
\end{equation}
for $k = 1,\dots,N_s$. The reduced basis $\Phibold$ is defined by compressing
the information in the snapshot matrix using POD, i.e., retaining the dominant
singular vectors of the translated snapshot matrix $\ubm - \ubm_0\onebold^T$
(to account for the affine offset).
This computationally expensive \textit{training phase} requires
$N_s$ solutions of the large-scale dynamical system and compression
of the resulting snapshot matrix of size $N_\ubm \times N_sN_t$, but
is only required once; the resulting reduced-order model can be
leveraged on a testing set $\Xibold^*$
without re-training to ameliorate the offline cost.
Generalizability of the basis to $\Xibold^*$ depends on the coverage of the
parameter space with training samples. Sophisticated greedy methods exist to
select training samples based on the maximum ROM error in the parameter space
\cite{rozza2008reduced, haasdonk2013convergence},
given a reliable error indicator is available.
Since we consider complex nonlinear problems, such error indicators with
high effectivity are not available so we simply use uniform
sampling---feasible in our setting due to low-dimensional parameter
spaces considered ($N_\mubold \leq 3$ in Section~\ref{sec:num-exp})---to
ensure sufficient coverage of the parameter space.

The computational cost of integrating the ROM (\ref{eqn:rom}) is dominated
by the evaluation of the reduced velocity function and its Jacobian and
the linear solve with the Jacobian matrix.
The computational complexity of evaluating the reduced velocity function
$(\taubold,t,\mubold)\mapsto\fbm_r(\taubold,t,\mubold)$ and its Jacobian
$(\taubold,t,\mubold)\mapsto\pder{\fbm_r}{\taubold}(\taubold,t,\mubold)$
are $\Ocal(g(N_\ubm)+N_\ubm k_\ubm)$ and
$\Ocal(g(N_\ubm)+N_\ubm k_\ubm^2)$, respectively.
An entire implicit time step requires
$\Ocal(g(N_\ubm)+N_\ubm k_\ubm^2+k_\ubm^3)$ operations, assuming
a direct solver is used for the linear system of equations,
which is almost exclusively the case due to the small size of
$k_\ubm$.
%Even though the computational cost of a ROM time step
%scales more favorably than that of the HDM, i.e., only linear dependence
%on the large dimension $N_\ubm$, the operations scaling with
%$N_\ubm$ constitute a significant bottleneck; substantial speedups
%will only be realized if the cost can be made independent of $N_\ubm$,
%the primary goal of hyperreduction methods (Section~\ref{sec:mor:hrom}).

%Complexity ROM: $f_r$ is $\Ocal(g(N_\ubm)+N_\ubm k_\ubm)$ and
%$df_r$ is $\Ocal(g(N_\ubm)+N_\ubm k_\ubm^2)$, where constant factors
%have been dropped and assume sparse Jacobian where number of nonzero entries
%per column does not scale with $N_\ubm$; entire implicit time step
%$\Ocal(g(N_\ubm)+N_\ubm k_\ubm^2 + k_\ubm^3)$

\subsection{Hyperreduction to accelerate projection of
            nonlinear terms}
\label{sec:mor:hrom}
Despite the substantial reduction in the dimensionality of the dynamical
system, from $N_\ubm$-dimensional in (\ref{eqn:semidisc}) to
$k_\ubm$-dimensional in (\ref{eqn:rom}) with $k_\ubm\ll N_\ubm$, it is
well-known the ROM only achieves marginal speedup relative to the HDM
due to an inherent bottleneck in the evaluation of the nonlinear term
$\fbm_r(\taubold, t, \mubold)$ with complexity proportional to the
large dimension $N_\ubm$: $\Ocal(g(N_\ubm)+N_\ubm k_\ubm)$. 
From the definition in (\ref{eqn:rom-velo}) it is clear that even though
$\fbm_r$ is a mapping between low-dimensional spaces, it is expensive to
evaluate due to a sequence of high-dimensional operations: reconstruction
of $\ubm_r = \bar{\ubm}+\Phibold\taubold$ ($\Ocal(N_\ubm k_\ubm)$
operations), evaluation of the HDM velocity function
$\fbm(\ubm_r, t, \mubold)$ ($\Ocal(g(N_\ubm))$ operations),
and projection of the velocity onto $\text{Range}(\Phibold)$
($\Ocal(N_\ubm k_\ubm)$ operations).
To overcome this computational bottleneck, a host of so-called
\textit{hyperreduction} methods
\cite{barrault2004empirical, ryckelynck2005priori, an2008optimizing,
      astrid2008missing, chaturantabut2010nonlinear, carlberg2011gnat,
      tiso2013discrete, farhat2015ecsw, zahr2017multilevel,
      yano2019discontinuous}
have been introduced to
approximate $\fbm_r$ at a cost that does not scale with $N_\ubm$.
However, these methods usually have limited prediction potential
for complex problems \cite{washabaugh2016phd} and, most importantly,
are difficult and code-intrusive to implement properly and achieve
substantial speedup in practice.

For example, the empirical interpolation method
\cite{barrault2004empirical} and its discrete variant
\cite{chaturantabut2010nonlinear} approximate the ROM velocity
function as
\begin{equation}
 \fbm_r(\taubold, t, \mubold) \approx
 \fbm_d(\taubold, t, \mubold) \coloneqq
 \Abm\Pbm^T\fbm(\bar{\ubm}+\Phibold\taubold, t, \mubold),
 \qquad 
  \Abm = \Psibold^T\Pibold(\Pbm^T\Pibold)^{-1}
         \in \Rbb^{k_\ubm\times k_\fbm},
\end{equation}
where $\Pibold \in \Rbb^{N_\ubm\times k_\fbm}$ is a basis for a
$k_\fbm$-dimensional subspace ($k_\fbm \ll N_\ubm$) used to approximate
the HDM velocity function $\fbm$ and
$\Pbm \in \Rbb^{N_\ubm \times k_\fbm}$ is a subset of the columns of
the $N_\ubm\times N_\ubm$ identity matrix used to sample entries
of the HDM velocity function. Then the (D)EIM reduced coordinates
$\func{\ybm_d}{[0, T]\times\Dcal}{\Rbb^{k_\ubm}}$ are computed such that
\begin{equation}
 \Mbm_r\oder{\ybm_d}{t} = \fbm_d(\ybm_d, t, \mubold), \qquad
 \ybm_d(0) = \ybm_0,
\end{equation}
and the HDM approximation
$\func{\ubm_d}{[0, T]\times\Dcal}{\Rbb^{N_\ubm}}$ is reconstructed as
\begin{equation} \label{eqn:deim-approx}
 \ubm(t, \mubold) \approx \ubm_d(t, \mubold) \coloneqq
 \bar\ubm + \Phibold\ybm_d(t, \mubold).
\end{equation}
As noted in
\cite{chaturantabut2010nonlinear, carlberg2013gnat, farhat2015ecsw},
for this approach to be efficient, it is not sufficient to first evaluate
$\fbm(\bar{\ubm}+\Phibold\taubold, t, \mubold)$ and then sample its entries.
Rather, the term $\Pbm^T\fbm$ must be evaluated directly given the appropriate
sampling of the state $\hat\Pbm^T\ubm$, where
$\hat\Pbm\in\Rbb^{N_\ubm \times k_s}$ is the matrix that samples all
entries of the state $\ubm$ required to evaluate the entries necessary
entries of the velocity function $\Pbm^T\fbm$. This approach assumes sparse
dependence of the velocity function on the state, i.e., each entry of the
velocity function depends on a small number of entries of the state vector.
This sparsity property is guaranteed if the dynamical system in
(\ref{eqn:semidisc}) corresponds to the semi-discretization of a PDE
using local methods e.g., finite difference or finite
volume methods. Direct implementation of $\Pbm^T\fbm$ given $\hat\Pbm^T\ubm$
in the context of a PDE discretization involves the use of a sparsified
or sample mesh on which the state $\hat\Pbm^T\ubm$ is defined
\cite{chaturantabut2010nonlinear}.
While successful, this approach is code intrusive and difficult to implement,
often requiring years of development to incorporate into existing codes.
In addition, the implementation is highly dependent on the semi-discretization
approach used for the PDE
\cite{chaturantabut2010nonlinear, carlberg2011gnat, tiso2013discrete,
      farhat2015ecsw, yano2019discontinuous}.
Other physics-based hyperreduction methods besides (D)EIM exist
\cite{barrault2004empirical, ryckelynck2005priori, an2008optimizing,
      astrid2008missing, chaturantabut2010nonlinear, carlberg2011gnat,
      tiso2013discrete, farhat2015ecsw, zahr2017multilevel,
      yano2019discontinuous}
to approximate the ROM velocity; however, they all rely on this concept of
partial assembly over a sample mesh and require a specialized, code-intrusive
implementation.

Assuming an efficient implementation of the sampled nonlinear velocity
function and its Jacobian, the computational complexity of evaluating
each term is $\Ocal(g(k_\fbm)+k_\fbm k_\ubm)$ and
$\Ocal(g(k_\fbm)+k_\fbm k_\ubm^2+\gamma k_\fbm k_\ubm)$
\cite{chaturantabut2010nonlinear}, respectively, where
$\gamma$ is the average number of nonzero entries per column of the
full Jacobian matrix $\pder{\fbm}{\xibold}$. In this work, we assume
$\gamma$ is independent of $N_\ubm$, which simplifies the complexity of the
reduced Jacobian evaluation to $\Ocal(g(k_\fbm)+k_\fbm k_\ubm^2)$. An
implicit time step requires $\Ocal(g(k_\fbm)+k_\fbm k_\ubm^2+k_\ubm^3)$
operations, assuming a direct solver is used for the linear system of
equations, which is independent of $N_\ubm$ and substantially cheaper than
integrating the HDM (\ref{eqn:semidisc}) and ROM without hyperreduction
(\ref{eqn:rom}).
%Complexity of DEIM: $f_r$ is $\Ocal(g(m)+m k_\ubm)$ and
%$df_r$ is $\Ocal(g(m)+m k_\ubm^2 + \gamma m k_\ubm)$, where constant factors
%have been dropped and assume sparse Jacobian where number of nonzero entries
%per column does not scale with $N_\ubm$; entire implicit time step
%$\Ocal(g(m)+m k_\ubm^2 + \gamma m k_\ubm + k_\ubm^3)$

\section{Non-intrusive hyperreduction using deep neural networks}
\label{sec:romnn}
We propose a new approach to hyperreduction that approximates the
ROM velocity function using a deep neural network, which we abbreviate
ROM-NN in the remainder. That is, we introduce a function
\begin{equation} \label{eqn:romnn-velo}
 \func{\hat\fbm_r}
      {\Rbb^{k_\ubm}\times[0, T]\times\Dcal\times\Rbb^{N_\wbm}}
      {\Rbb^{k_\ubm}}, \qquad
 (\taubold, t, \mubold, \nubold) \mapsto
 \hat\fbm_r(\taubold, t, \mubold; \nubold)
\end{equation}
and vector of weights $\wbm \in \Rbb^{N_\wbm}$ such that
\begin{equation} \label{eqn:romnn-velo-approx}
 \hat\fbm_r(\ybm(t, \mubold), t, \mubold; \wbm) \approx
 \fbm_r(\ybm(t, \mubold), t, \mubold)
\end{equation}
for any $t \in [0, T]$ and $\mubold \in \Dcal$ and
compute $\func{\ybm_n}{[0, T]\times\Dcal}{\Rbb^{k_\ubm}}$
that solves
\begin{equation}
    \Mbm_r\oder{\ybm_n}{t} = \hat{\fbm}_r(\ybm_n,t,\mubold), \qquad
    \ybm_n(0) = \ybm_0.
\end{equation}
The HDM approximation
$\func{\ubm_n}{[0, T]\times\Dcal}{\Rbb^{N_\ubm}}$ is reconstructed as
\begin{equation} \label{eqn:romnn-approx}
 \ubm(t, \mubold) \approx \ubm_n(t, \mubold) \coloneqq
 \bar\ubm + \Phibold\ybm_n(t, \mubold).
\end{equation}
If (\ref{eqn:romnn-velo-approx}) holds, we expect $\ybm_n$
to be a good approximation to $\ybm$ and, provided the reduced basis is
sufficient, $\ubm_n$ to be a good approximation to $\ubm$.

The ROM velocity (\ref{eqn:rom-velo}) is a $k_\ubm$-valued function of
$k_\ubm+1+N_\mubold$ variables; training a neural network to
approximate this mapping requires a (large) number of instances
of the function input
$(\taubold, t, \mubold) \in \Rbb^{k_\ubm}\times[0, T]\times\Dcal$
and the corresponding output $\fbm_r(\taubold, t, \mubold) \in\Rbb^{k_\ubm}$
so the network weights can be tuned to minimize a loss function.
In this work, the network weights are defined as the solution of the
following optimization problem
\begin{equation} \label{eqn:opt-loss}
 \optunc{\wbm\in\Rbb^{N_\wbm}}{\frac{1}{2}
               \sum_{i=1}^{N_t}
               \sum_{j=1}^{N_s}
               \norm{\hat\fbm_r(\taubold_{ij}, t_i, \mubold_j; \wbm) -
                    \fbm_r(\taubold_{ij}, t_i, \mubold_j)}_2^2},
\end{equation}
where $\{t_i\}_{i=1}^{N_t} \subset [0, T]$ are
the nodes of the temporal discretization (Section~\ref{sec:mor:hdm}),
$\{\mubold_k\}_{k=1}^{N_s} \subset \Dcal$ are the training parameters
(Section~\ref{sec:mor:rom}), and $\taubold_{ij} \in \Rbb^{k_\ubm}$
for $i = 1,\dots,N_t$, $j = 1, \dots, N_s$ are the reduced coordinates
used for training the network.
Given the requirement in (\ref{eqn:romnn-velo-approx}) that the DNN velocity
function matches the ROM velocity function \textit{on the manifold of ROM
solutions}, a sensible choice is $\taubold_{ij} = \ybm(t_i,\mubold_j)$.
While consistent with the requirement in (\ref{eqn:romnn-velo-approx}),
this approach requires that both the HDM solution $\ubm(t, \mubold)$ and
expensive ROM solution $\ybm(t, \mubold)$, i.e., without hyperreduction,
be computed for each $\mubold \in \Xibold_0$ to define the training data,
which can substantially increase the offline cost. To mitigate the additional
cost of computing the expensive ROM solution, we propose to use the projection
of the HDM solution onto the subspace $\Vcal$ in place of the ROM solution
itself. That is, we take
$\taubold_{ij} = \tilde\ybm(t_i, \mubold_j)$, where
\begin{equation}
 \func{\tilde\ybm}{[0, T]\times\Dcal}{\Rbb^{k_\ubm}}, \qquad
 (t, \mubold) \mapsto \Phibold^T(\ubm(t, \mubold)-\bar\ubm),
\end{equation}
which requires exactly the \textit{same} data used to compute the reduced
basis $\Phibold$. The complete training procedure is summarized in
Algorithm~\ref{alg:romnn-offline}.
\begin{algorithm}
 \caption{Training procedure for deep learning-based reduced-order model}
 \label{alg:romnn-offline}
 \begin{algorithmic}[1]
  %\REQUIRE Training set $\Xibold_0 \subset \Dcal$ ($|\Xibold_0|= N_s$),
  \REQUIRE Training set $\{\mubold_j\}_{j=1}^{N_s} \subset \Dcal$,
           temporal discretization $\{t_i\}_{i=1}^{N_t} \subset [0, T]$,
           initial condition $\ubm_0 \in \Rbb^{N_\ubm}$
  \ENSURE Affine offset $\bar\ubm \in \Rbb^{N_\ubm}$,
          reduced basis $\Phibold\in\Rbb^{N_\ubm\times k_\ubm}$,
          network weights $\wbm \in \Rbb^{N_\wbm}$ 
  \FOR{$j = 1, \dots, N_s$}
    \STATE Compute solution of the HDM dynamical system
           $\ubm(\,\cdot\,,\,\mubold_j)$
  \ENDFOR
  \STATE Define the snapshot matrix $\Xbm\in\Rbb^{N_\ubm\times N_sN_t}$
         according to  (\ref{eqn:snapmat})
  \STATE Compute the left singular vectors of $\Xbm - u_0\onebold^T$:
         $\ubm_i \in \Rbb^{N_\ubm}$, $i = 1,\dots,N_sN_t$
  \STATE Define reduced subspace:
  $\bar\ubm \leftarrow \ubm_0$,
  $\Phibold \leftarrow \begin{bmatrix} \ubm_1 & \cdots & \ubm_r \end{bmatrix}$
  \FOR{$j = 1, \dots, N_s$}
    \FOR{$i = 1, \dots, N_t$}
      \STATE Compute $\taubold_{ij} = \Phibold^T(\ubm(t_i,\mubold_j)-\bar\ubm)$
    \ENDFOR
  \ENDFOR
  \STATE Solve (\ref{eqn:opt-loss}) for network weights $\wbm$
 \end{algorithmic}
\end{algorithm}

\begin{figure}
 \centering
 \raisebox{-0.5\height}{
   \includegraphics[width=0.3\textwidth]{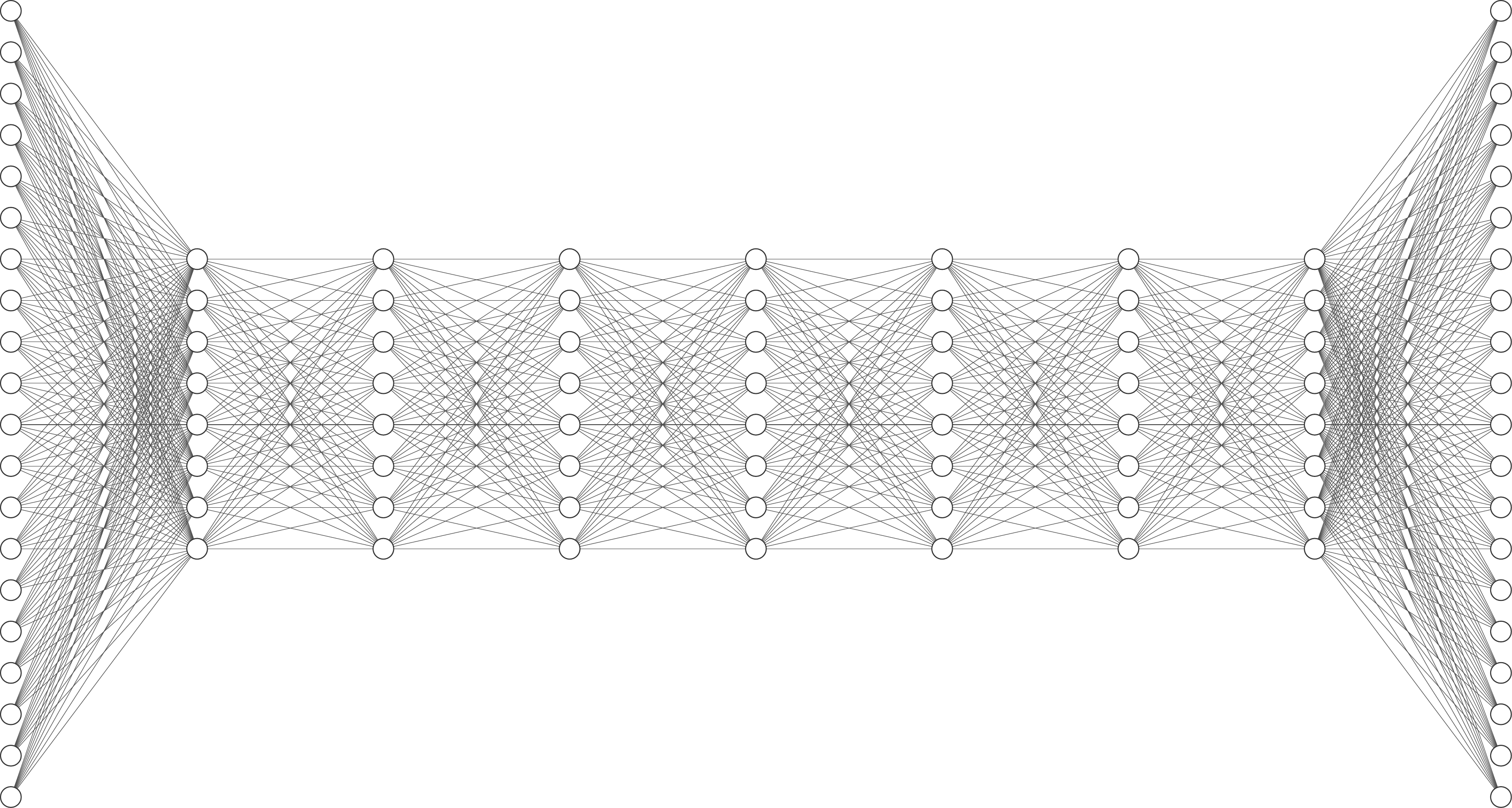}
 } \quad
 \raisebox{-0.5\height}{
   \includegraphics[width=0.3\textwidth]{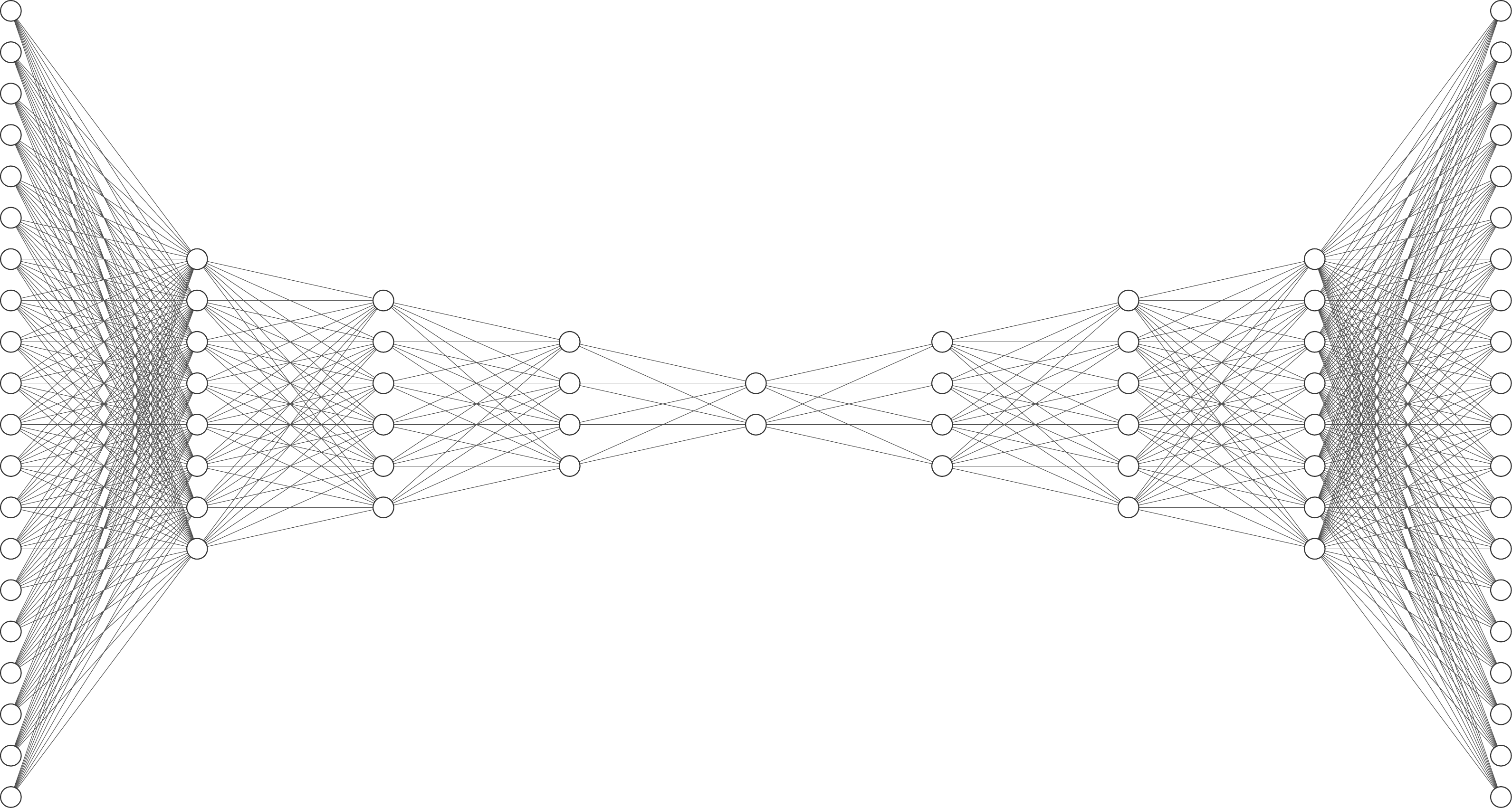}
 } \quad
 \raisebox{-0.5\height}{
   \includegraphics[width=0.3\textwidth]{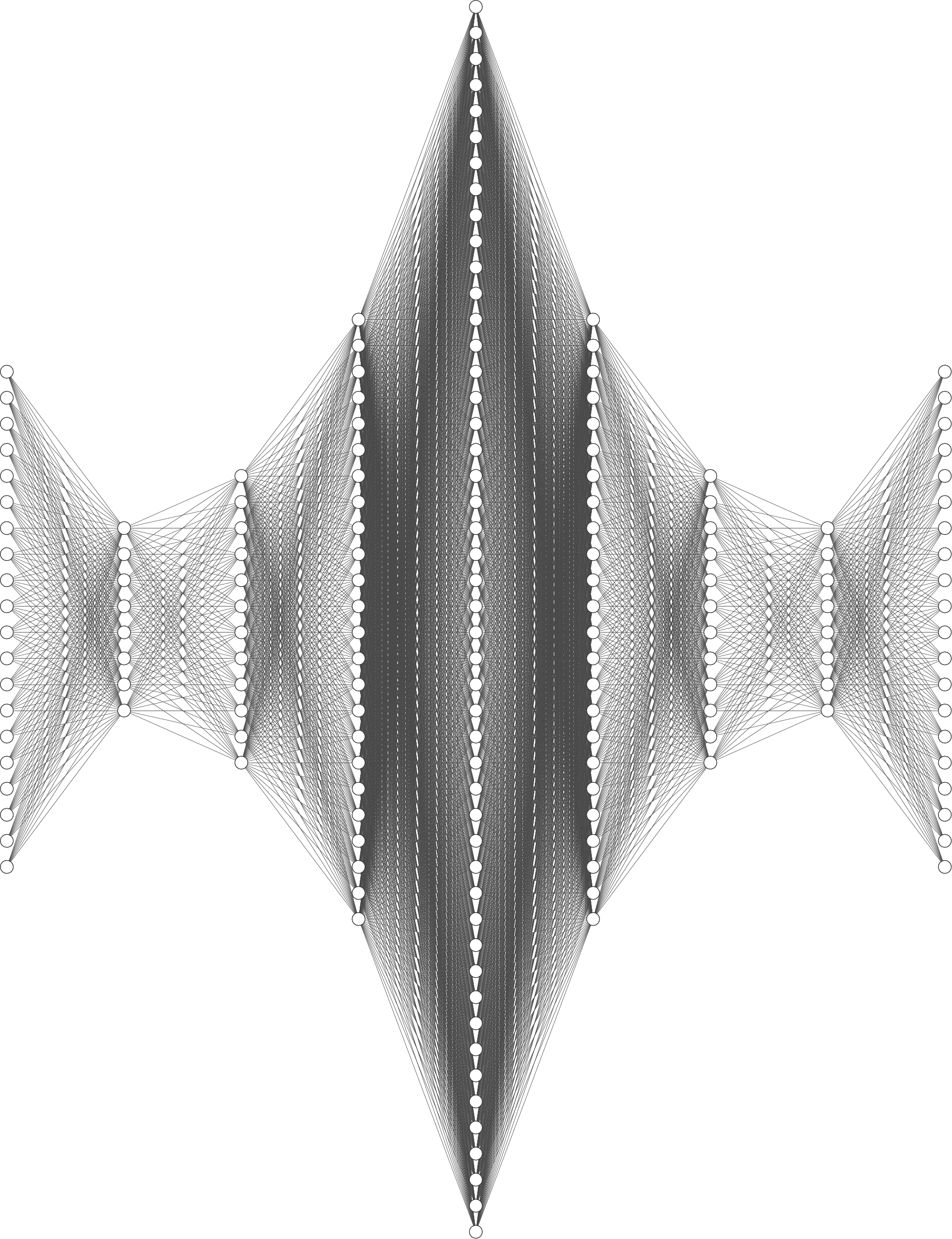}
 }
 \caption{Neural network structures tested: uniform (\textit{left}),
          converging-diverging (\textit{middle}), and diverging-converging
          (\textit{right}).}
 \label{fig:nnstruct}
\end{figure}

In this work, we define $\hat\fbm_r$ using a fully-connected, feed-forward
neural network (FCNN) architecture, which contains one input layer (the input
vector $(\taubold, t, \mubold)$), five hidden layers, and one output layer
(the prediction $\hat\fbm_r(\taubold,t,\mubold)$). % (Figure~\ref{fig:nn}).
Each layer is fed forward to the next layer by a linear weighted sum and
nonlinear activation function (e.g., reLU). The network is built in the form
of sparse autoencoder (SAE), namely the hidden layers follow a decoder-encoder
structure in order to capture the complex hidden nonlinear pattern of the
mapping. The number of neurons for each layers from the input to the output are 
$(80, 120, 240, 480, 240, 120, 80)$. Standardized normalization is applied for
both input and output layers. The training is conducted in a supervised manner,
i.e., minimizing the loss function  of data misfit, using a stochastic
gradient descent (SGD) based optimizer (e.g. Adam
algorithm~\cite{kingma2014adam}).  To avoid over-fitting, the
dropout~\cite{srivastava2014dropout} and early stopping~\cite{yao2007early}
techniques are applied.  To demonstrate the robustness of the ROM-NN, the
architecture and hyperparameters of the network remain the same for all test
cases throughout the paper.
We considered a number of other FCNN structures (uniform,
 converging-diverging, diverging-converging; Figure~\ref{fig:nnstruct})
 of varying depths and breadths and found the accuracy of the ROM-NN
 to be rather insensitive to the network structure; for the problems
 considered, even a shallow neural network with only two hidden layers
 and $80$ neurons per layer leads to a ROM-NN with similar overall accuracy
 as one with the aforementioned diverging-converging structure. The only
 exception being in the limit of sparse training where a deeper network leads
 to a slightly more accurate ROM-NN.
%Note that a more comprehensive study on optimizing the neural network structure
%is out of the scope of this work.

This approach is guaranteed to be \textit{non-intrusive} because the
training procedure only relies on snapshots of the HDM solutions and
evaluations of the ROM velocity function and the online solution only
requires evaluation of the neural network velocity function (forward
pass) and its derivative with respect to $\taubold$ (backward
propagation). As a result, both the DNN and dynamical system code
can be treated as \textit{black boxes}, which substantially eases
the implementation burden. Another advantage of the proposed ROM-NN
method is we directly approximate $\fbm_r$, a mapping between
low-dimensional input and output spaces, using nonlinear basis functions.
We will show in our numerical experiments (Section~\ref{sec:num-exp}) that
this approximation, when sufficiently trained, mitigates some parametric
robustness and stability issues of traditional hyperreduction techniques,
such as (D)EIM, that approximate the mapping
$(\taubold, t, \mubold) \mapsto \fbm(\bar\ubm+\Phibold\taubold, t, \mubold)$
(low-dimensional input space, high-dimensional output space) using a
linear basis.

The computational complexity of a single pass through a FCNN with
$M+2$ layers with the $i$th layer consisting of $m_i$ neurons,
$i = 0,\dots,M+1$, where layer $0$ is the input and layer $M+1$ is the
output, is $\Ocal(\sum_{i=1}^{M+1} m_{i-1} m_i)$,
which can be seen from a simple analogy to dense matrix-vector
multiplication. In our case, $m_0 = k_\ubm+N_\mubold+1$
and $m_{M+1} = k_\ubm$. Therefore evaluation of the approximate velocity
function
$(\taubold,t,\mubold,\nubold)\mapsto\hat{f}_r(\taubold,t,\mubold,\nubold)$
scales quadratically with the breadth of the network and linearly in the depth:
$\Ocal(k_\ubm m_1 + N_\mubold m_1 + k_\ubm m_M + \sum_{i=2}^M m_{i-1} m_i)$.
In the special case where all hidden layers are the same size
$m_i = m$ for $i = 1, \dots, M$, this reduces to
$\Ocal(k_\ubm m + N_\mubold m + Mm^2)$. The cost to evaluate the Jacobian
$(\taubold,t,\mubold,\nubold)\mapsto\pder{\hat{f}_r}{\taubold}(\taubold,t,\mubold,\nubold)$
using automatic differentiation is within a constant factor of the cost for
the velocity function itself \cite{baur1983complexity, griewank2008evaluating}.
Therefore the cost of an entire implicit time step is
$\Ocal(k_\ubm^3 + k_\ubm m_1 + N_\mubold m_1 + k_\ubm m_M + \sum_{i=2}^M m_{i-1} m_i)$ assuming a direct solver is used for the linear system. In the special
case where all hidden layers are the same size, this reduces to
$\Ocal(k_\ubm^3+k_\ubm m + N_\mubold m + Mm^2)$. Assuming the breadth
($m$) and depth ($M$) of the network scale independently of $N_\ubm$,
we expect this approach to be efficient because the dependence on the
large dimension has been removed. If we require the network breadth and
depth to be on the order of the size of the reduced basis, i.e.,
$m \sim \Ocal(k_\ubm)$ and $M \sim \Ocal(k_\ubm)$, the complexity of a time
step reduces to $\Ocal(k_\ubm^3+N_\mubold k_\ubm)$. This shows the bottleneck
caused by the nonlinear terms has been completely eliminated using the
proposed FCNN approximation of the nonlinear velocity function with
$\Ocal(k_\ubm)$ layers and neurons per layer because, asymptotically, the cost
is similar to that of a direct solve with the reduced Jacobian matrix. A
similar result follows with broader, shallower networks, e.g.,
$m \sim \Ocal(k_\ubm^{3/2})$ and $M \sim \Ocal(1)$, with complexity
per time step: $\Ocal(k_\ubm^3+N_\mubold k_\ubm^{3/2})$. Similarly,
deeper, narrower networks can be used, e.g.,
$m \sim \Ocal(k_\ubm^{1/2})$ and $M \sim \Ocal(k_\ubm^2)$,
with a complexity per time step of $\Ocal(k_\ubm^3+N_\mubold k_\ubm^{1/2})$.
The number of parameters ($N_\mubold$) is usually small and always
independent of $k_\ubm$, therefore the dominant complexity for all the
network structures mentioned is $\Ocal(k_\ubm^3)$.

%The
%reduced cost and non-intrusivity of the method come replacing queries to the
%original nonlinear velocity function with FCNN queries

%Complexity of ROM-NN: consider FCNN with $d+1$ levels with each
%level containing $N_i$ nodes $i = 0, \dots, d$, where level $0$
%is the input layer and level $d$ is the output layer.
%$N_0 = k_\ubm+N_\mubold$ adn $N_{d} = k_\ubm$.
%Owning to connection between FCNN and matrix multiplication, computational
%cost of forward prop through newtork is
%$\Ocal((k_\ubm+N_\ubm) N_1 + N_{d-1}k_\ubm + \sum_{i=1}^{d-1} N_iN_{i+1})$.
%If all hidden layers have same number of neurons:
%$\Ocal(k_\ubm N + N_\ubm N + d N^2)$.
%Therefore residual/Jacobian evaluations
%cost is $\Ocal(\sum_{i=0}^d N_iN_{i+1})$ and entire implicit step is
%$\Ocal(k_\ubm^3+\sum_{i=0}^d N_iN_{i+1})$.

%Besides the non-intrusivity of these approach, an advantage over
%(D)EIM is we directly approximate $\fbm_r$ which is a mapping from
%a low-dimensional input and output space using nonlinear basis functions,
%whereas (D)EIM directly approximate $\fbm()$, which is a low-dimensional
%input space but a high-dimensional output space using linear basis functions.
%
%\begin{itemize}
% \item Since we are approximating a small input-output
%mapping ($k_\ubm+1+N_\mubold$ inputs, $k_\ubm$ outputs) using the DNN, we
%expect the results to be accurate beyond the training set. Furthermore, we
%expect this to be superior to traditional hyperreduction techniques that
%directly approximate the high-dimensional nonlinear term
%$(\ybm,\mubold) \mapsto \fbm(\Vbm\ybm,\mubold)$.
%\end{itemize}

\section{Numerical experiments}
\label{sec:num-exp}
In this section, we test the accuracy, stability, and parametric robustness
of the proposed ROM-NN method using two dynamical systems that result from
the semi-discretization of nonlinear, hyperbolic PDEs. We compare the
performance of the ROM-NN method to a standard Galerkin-POD ROM, which
provides a theoretical lower bound on the accuracy of the ROM-NN, and
the most popular intrusive hyperreduction method, (D)EIM. For both
problems, we define the parameter space $\Dcal$ and introduce two
subsets $\Xibold_0 \subset \Xibold^* \subset \Dcal$, where $\Xibold_0$
are the parameters used to train the reduced-order models and
$\Xibold^*$ are all parameters where the accuracy of the models
is tested (includes the training points).  Recall the definition
of the parametric HDM solution $\ubm$ (\ref{eqn:semidisc}) and its
approximation provided by the ROM $\ubm_r$ (\ref{eqn:rom-ansatz}),
(D)EIM $\ubm_d$ (\ref{eqn:deim-approx}), and ROM-NN $\ubm_n$
(\ref{eqn:romnn-approx}).
For a given parameter $\mubold \in \Dcal$, we
quantify the error between the HDM solution $\ubm(\,\cdot\,,\,\mubold)$
and an approximate solution $\vbm(\,\cdot\,,\,\mubold)$ as
\begin{equation}
 \epsilon(\vbm; \mubold) \coloneqq
 \sqrt{\frac{\sum_{i=1}^{N_t} \norm{\vbm(t_i, \mubold)-\ubm(t_i, \mubold)}^2}
            {\sum_{i=1}^{N_t} \norm{\ubm(t_i, \mubold)}^2}}.
\end{equation}
Therefore the error in the ROM, (D)EIM, and ROM-NN solutions are
\begin{equation}
 \epsilon_r(\mubold) \coloneqq \epsilon(\ubm_r; \mubold), \qquad
 \epsilon_d(\mubold) \coloneqq \epsilon(\ubm_d; \mubold), \qquad
 \epsilon_n(\mubold) \coloneqq \epsilon(\ubm_n; \mubold),
\end{equation}
respectively. In the rest of this section, we will consider the
statistics (minimum, maximum, and median) of these error metrics
over the training set $\Xibold_0$ and testing set
$\Xibold^*\setminus\Xibold_0$.
In this work, we do not compare the computational cost of the HDM, ROM,
(D)EIM, and ROM-NN because it is heavily dependent on the implementation
and a number algorithmic choices, e.g., choice of linear solver.

\subsection{One-dimensional viscous Burgers' equation}
The first numerical experiment we consider is solution of the
one-dimensional, parametrized, viscous Burgers' equation in the
domain $\Omega \coloneqq (0, 1)$, where
$\func{u}{\Omega \times [0, T]\times\Dcal}{\Rbb}$ solves
\label{sec:numexp:vburg}
\begin{equation}
 \begin{aligned}
  \partial_t u(x, t, \mubold) +
  u(x, t, \mubold)\partial_x u(x, t, \mubold) &=
  \nu(\mubold) \partial_{xx} u(x, t, \mubold) + g(x, \mubold),
  &&x \in (0, 1), &&t \in [0, T], &&\mubold \in \Dcal, \\
  u(0, t, \mubold) &= 0, \quad
  u(1, t, \mubold)  = 0,
  && &&t \in [0, T], && \mubold \in \Dcal.
 \end{aligned}
\end{equation}
The time interval is taken as $T = 1$ and the parameter space is
$\Dcal \coloneqq [0.01, 0.1] \times [2, 3] \times [0, 1] \subset \Rbb^3$.
For any $\mubold \in \Dcal$ where
$\mubold = (\mu_1, \mu_2, \mu_3)$,
the parametrized viscosity and source term are defined as
\begin{equation}
 \nu(\mubold) = \mu_1, \quad
 g(x, \mubold) = \mu_2 e^{\mu_3 x}.
\end{equation}
The PDE is discretized in space using $200$ linear finite elements
with essential boundary conditions strongly enforced to yield a
dynamical system of the form (\ref{eqn:semidisc}) with a total of
$N_\ubm = 199$ spatial degrees of freedom. The dynamical system
is discretized in time using the two-stage diagonally implicit
Runge-Kutta method \cite{alexander1977diagonally} with $100$ time steps.

For this problem, we define the testing set $\Xibold^*$ as the uniform
sampling of $\Dcal$ on a $5 \times 5 \times 5$ grid for a total of
$|\Xibold^*| = 125$ parameter configurations. We consider two training sets:
$\Xibold_0^a$, $\Xibold_0^b$ are the uniform samplings of $\Dcal$
on a $2\times 2\times 2$ and $3\times 3\times 3$ grid, respectively.
By construction, $\Xibold^a \subset \Xibold^b \subset \Xibold^*$.
For both training sets, we construct a POD-Galerkin ROM without
hypereduction, accelerated with (D)EIM, and accelerated with
the neural network approximation of $\fbm_r$ using $k_\ubm = 8$
and test each model on all points in $\Xibold^*$.

The reduced-order model without hyperreduction is the most stable and
accurate method, which is expected since it computes the velocity
function $\fbm_r$ exactly. Nonlinear approximation via (D)EIM is the
least accurate approach and even goes unstable for a number of training
and testing points for this small basis size. The neural network approach to
approximate the nonlinear terms is more accurate than (D)EIM and is
stable for all points in $\Xibold^*$. These observations are
taken from Figures~\ref{fig:vburg_train0_stvc}~and~\ref{fig:vburg_train1_stvc},
which contain the PDE state vector computed with each model at several
instances in time for various points in $\Xibold^*$ and
Figure~\ref{fig:vburg_train_all_err0} that directly compares the
accuracy of (D)EIM and ROM-NN.

\ifbool{fastcompile}{}{
\begin{figure}
 \centering
 \input{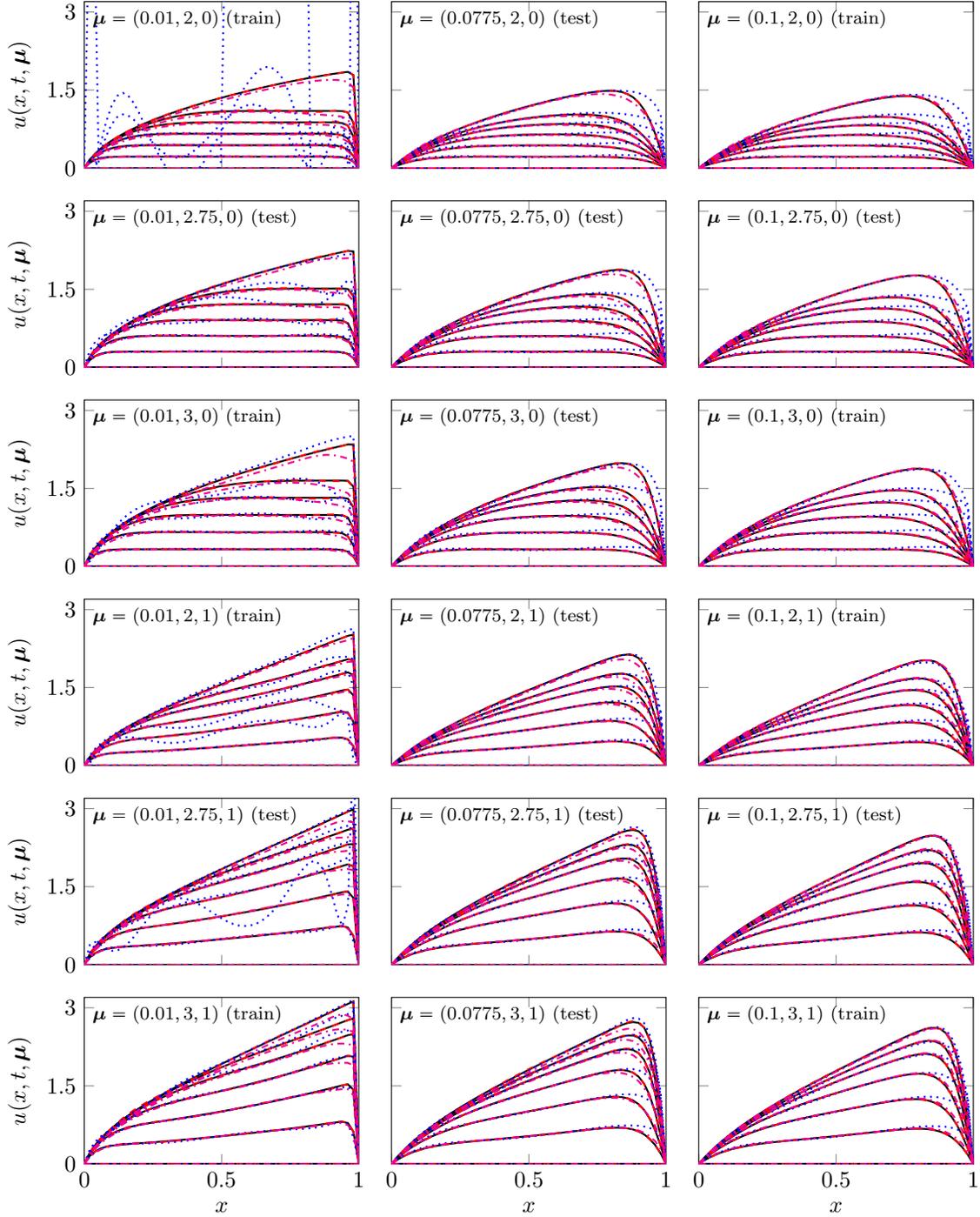}
 \caption{Snapshots of viscous Burgers' equation
          ($t = 0, 0.11, 0.22, 0.33, 0.44, 0.55, 1$)
          at various parameter configurations using HDM (\ref{line:hdm_stvc}),
          ROM (\ref{line:rom_stvc}), (D)EIM (\ref{line:deim_stvc}),
          ROM-NN (\ref{line:romnn_stvc}). The model reduction methods
          are trained using $8$ parameter samples ($\Xibold_0^a$)
          for a total of $800$ snapshots and compressed to a size
          $k_\ubm = 8$. In most cases, including both training and testing
          configuration, the ROM-NN model is more accurate than the (D)EIM
          model and does not exhibit the same stability issues.}
 \label{fig:vburg_train0_stvc}
\end{figure}
}

\ifbool{fastcompile}{}{
\begin{figure}
 \centering
 \input{tikz/vburg_train1_stvc.tikz}
 \caption{Snapshots of viscous Burgers' equation
          ($t = 0, 0.11, 0.22, 0.33, 0.44, 0.55, 1$)
          at various parameter configurations using HDM (\ref{line:hdm_stvc}),
          ROM (\ref{line:rom_stvc}), (D)EIM (\ref{line:deim_stvc}),
          ROM-NN (\ref{line:romnn_stvc}). The model reduction methods
          are trained using $27$ parameter samples ($\Xibold_0^b$)
          for a total of $2700$ snapshots and compressed to a size
          $k_\ubm = 8$. In most cases, including both training and testing
          configuration, the ROM-NN model is more accurate than the (D)EIM
          model and does not exhibit the same stability issues.}
 \label{fig:vburg_train1_stvc}
\end{figure}
}

For the training set $\Xibold_0^a$, the minimum error across both
the training and testing sets are comparable for (D)EIM and
ROM-NN. Since (D)EIM is unstable on both training and testing
points, the maximum error is large. The ROM-NN approach is stable
for all points in $\Xibold^*$; however, its maximum error on the
test set $\Xibold^*\setminus\Xibold_0^a$ is large ($\approx 22\%$).
The median error for the ROM-NN less than $3\%$ on the training
set and $5\%$ on the testing set, while (D)EIM has median errors
up to five times larger ($16\%$ training, $11\%$ testing)
(Table~\ref{tab:vburg_train0_err}).
By increasing the training set to $\Xibold_0^b$ and keeping the ROM size
fixed ($k_\ubm = 8$), the stability of (D)EIM further degrades (15
unstable points in $\Xibold^*$), but the accuracy improves for stable
configurations (median errors decrease). The ROM-NN training errors are
similar to the case where $\Xibold_0^a$ is used as the training set,
but the errors on the testing set become smaller suggesting the additional
training leads to better prediction. The median and maximum errors of
the ROM-NN for both training and testing sets are roughly $3\%$ and
$6\%$, respectively (Table~\ref{tab:vburg_train1_err}).

\begin{table}
 \centering 
 \caption{Summary of the performance of the model reduction methods
          trained on $\Xibold_0^a$, compressed to $k_\ubm = 8$,
          and tested on $\Xibold^*$. The error statistics are reported
          for the training set $\Xibold_0^a$ and testing set
          $\Xibold^*\setminus\Xibold_0^a$ separately.
          The ROM-NN is stable for all training and testing points considered
          and has a median error less than $3\%$ on the training set
          and about $5\%$ on the testing set. The (D)EIM method goes
          unstable at a number of training and testing point and has a median
          error greater than $10\%$.}
 \label{tab:vburg_train0_err}
 \begin{tabular}{r|cc|cc|cc}
  & \multicolumn{2}{c|}{ROM} & \multicolumn{2}{c|}{(D)EIM} & \multicolumn{2}{c}{ROM-NN} \\\hline
  & Train set & Test set & Train set & Test set & Train set & Test set \\\hline
  Unstable ($\#$) & 0 & 0 & 1 & 5 & 0 & 0 \\
  Minimum error & 2.58e-03 & 1.48e-03 & 2.81e-02 & 3.07e-02 & 1.68e-02 & 2.28e-02 \\
  Maximum error & 7.05e-03 & 2.54e-02 & 1.99e+01 & 5.49e+02 & 6.17e-02 & 2.25e-01 \\
  Median error & 4.66e-03 & 5.84e-03 & 1.58e-01 & 1.12e-01 & 2.71e-02 & 5.24e-02 \\
\end{tabular}

\end{table}

\begin{table}
 \centering 
 \caption{Summary of the performance of the model reduction methods
          trained on $\Xibold_0^b$, compressed to $k_\ubm = 8$,
          and tested on $\Xibold^*$. The error statistics are reported
          for the training set $\Xibold_0^b$ and testing set
          $\Xibold^*\setminus\Xibold_0^b$ separately.
          The ROM-NN is stable for all training and testing points considered
          and has a median error less than $3\%$ on both the training
          and testing set. The (D)EIM method goes unstable at a number of
          training and testing point and has a median error greater than
          $10\%$.}
 \label{tab:vburg_train1_err}
 \begin{tabular}{r|cc|cc|cc}
  & \multicolumn{2}{c|}{ROM} & \multicolumn{2}{c|}{(D)EIM} & \multicolumn{2}{c}{ROM-NN} \\\hline
  & Train set & Test set & Train set & Test set & Train set & Test set \\\hline
  Unstable ($\#$) & 0 & 0 & 6 & 9 & 0 & 0 \\
  Minimum error & 1.43e-03 & 1.54e-03 & 2.95e-02 & 9.39e-03 & 8.09e-03 & 7.58e-03 \\
  Maximum error & 8.62e-03 & 1.81e-02 & 2.96e+02 & 1.65e+03 & 6.26e-02 & 5.99e-02 \\
  Median error & 2.67e-03 & 2.51e-03 & 1.22e-01 & 1.16e-01 & 2.78e-02 & 2.47e-02 \\
\end{tabular}

\end{table}

\ifbool{fastcompile}{}{
\begin{figure}
 \centering
 \begin{tikzpicture}
\begin{groupplot} [
group style={group size = 2 by 1}]
\nextgroupplot[ymax=1.0, xmode=log, height=0.35\textwidth, width=0.45\textwidth, xtick={1e-2, 1e-1, 1e0}, ytick={1e-2, 1e-1, 1e0}, xlabel={$\epsilon_d(\mubold)$}, xmax=1.0, ylabel={$\epsilon_n(\mubold)$}, xmin=0.01, ymode=log, ymin=0.01]
\addplot [thick, black, solid]
coordinates {
(  0.01000000,   0.01000000)
(  1.00000000,   1.00000000)};\label{line:eye}

\addplot [black, only marks, mark=x, mark size=1.25]
coordinates {
( 19.89929500,   0.05500041)
(  0.19307856,   0.01675903)
(  0.21899561,   0.02852714)
(  0.09164970,   0.02390333)
(  0.13513397,   0.06172398)
(  0.02814268,   0.04156330)
(  0.18035288,   0.02448067)
(  0.08029719,   0.02577013)};\label{line:deim_v_romnn_train}

\addplot [red, only marks, mark=o, mark size=1]
coordinates {
(  7.53724444,   0.03243887)
(  0.04451748,   0.15207857)
( 548.59069655,   0.22541810)
(  0.08480804,   0.08498527)
(  0.06419097,   0.06162389)
(  0.07398385,   0.09895198)
(  0.12606755,   0.10803726)
(  0.21722182,   0.06658168)
(  0.15170692,   0.06836829)
(  0.12600257,   0.05391446)
(  0.10183911,   0.06216757)
(  0.07894473,   0.05601800)
(  0.05850542,   0.06607434)
(  0.19127714,   0.03571803)
(  0.15957370,   0.03834809)
(  0.13211785,   0.04676576)
(  0.10866151,   0.03674300)
(  0.08843674,   0.03615380)
(  0.17715809,   0.04449585)
(  0.14161362,   0.05908679)
(  0.11323283,   0.05178211)
( 387.57991409,   0.02727517)
( 490.66797819,   0.04593684)
(  1.01032089,   0.18950764)
(  0.06741370,   0.20295049)
(  0.04348321,   0.02742724)
(  0.07337130,   0.06894298)
(  0.06814337,   0.06908558)
(  0.10701867,   0.07631169)
(  0.18298270,   0.09657746)
(  0.36665822,   0.07429761)
(  0.14029597,   0.06509119)
(  0.11509856,   0.05835839)
(  0.09113602,   0.05680340)
(  0.06872741,   0.05815122)
(  0.05182895,   0.06866247)
(  0.18020759,   0.03918855)
(  0.15028403,   0.04159539)
(  0.12429127,   0.03620335)
(  0.10175478,   0.03392908)
(  0.08184842,   0.03626420)
(  0.20740396,   0.02910227)
(  0.16830914,   0.04163240)
(  0.13548082,   0.03955048)
(  0.10926060,   0.03740977)
(  0.08892738,   0.02277587)
(  0.10601826,   0.02562059)
(  0.11727648,   0.04769588)
(  0.03340412,   0.17207366)
(  0.03818849,   0.10300901)
( 67.59511661,   0.02976386)
(  0.06881975,   0.06764146)
(  0.08890143,   0.06359034)
(  0.15399817,   0.07675333)
(  0.20657262,   0.07404292)
(  0.58748044,   0.07368601)
(  0.12994295,   0.06491183)
(  0.10505261,   0.06074733)
(  0.08131589,   0.05681791)
(  0.06038034,   0.05735640)
(  0.05089135,   0.06820237)
(  0.17041417,   0.04068006)
(  0.14198058,   0.03977138)
(  0.11709860,   0.03432369)
(  0.09515952,   0.03175229)
(  0.07539026,   0.03489719)
(  0.19725208,   0.03296883)
(  0.16063917,   0.03569899)
(  0.13009942,   0.03737311)
(  0.10558308,   0.04244625)
(  0.08617403,   0.02419497)
(  0.26243485,   0.03460209)
(  0.11152925,   0.02765833)
(  0.51850822,   0.07953060)
(  0.05782379,   0.04425536)
(  0.08391899,   0.03316404)
(  0.07648307,   0.07313379)
(  0.12291812,   0.06392892)
(  0.36841491,   0.06431689)
(  0.25587000,   0.06672647)
( 466.12616388,   0.07061361)
(  0.12041838,   0.06776267)
(  0.09573108,   0.06263567)
(  0.07253368,   0.05463569)
(  0.05498067,   0.05601377)
(  0.05668777,   0.06610793)
(  0.16165233,   0.03965750)
(  0.13443676,   0.03628212)
(  0.11037693,   0.03404298)
(  0.08880827,   0.03205717)
(  0.06913135,   0.03474061)
(  0.18830405,   0.03251421)
(  0.15389034,   0.02855406)
(  0.12525619,   0.04960850)
(  0.10207051,   0.04855550)
(  0.08330771,   0.02498474)
(  0.03327049,   0.02738485)
(  0.04208729,   0.03180018)
(  0.03065181,   0.02393556)
(  0.09665417,   0.08431355)
(  0.73900435,   0.06946895)
(  0.19911320,   0.06928054)
(  0.31729343,   0.06598272)
(  0.81353293,   0.06902922)
(  0.11156487,   0.07033519)
(  0.08710634,   0.05653851)
(  0.06517745,   0.05238665)
(  0.05381967,   0.05552969)
(  0.06809087,   0.06793570)
(  0.15372123,   0.03750470)
(  0.12748343,   0.03082606)
(  0.10401811,   0.03265645)
(  0.08267801,   0.02967007)
(  0.06319144,   0.03823515)
(  0.14786106,   0.03976691)
(  0.12080095,   0.05482355)
(  0.09864049,   0.04556620)};\label{line:deim_v_romnn_test}

\nextgroupplot[yticklabels={,,}, ymax=1.0, xmode=log, height=0.35\textwidth, width=0.45\textwidth, xtick={1e-2, 1e-1, 1e0}, ytick={1e-2, 1e-1, 1e0}, xlabel={$\epsilon_d(\mubold)$}, xmax=1.0, xmin=0.01, ymode=log, ymin=0.01]
\addplot [thick, black, solid]
coordinates {
(  0.01000000,   0.01000000)
(  1.00000000,   1.00000000)};\label{line:eye}

\addplot [black, only marks, mark=x, mark size=1.25]
coordinates {
(  0.29528309,   0.02244362)
( 37.35409663,   0.02098433)
( 80.97647587,   0.02276096)
(  0.02953009,   0.02241784)
(  0.04048908,   0.01417453)
(  0.09150710,   0.02775870)
(  0.18486872,   0.01422387)
(  0.03719423,   0.03334434)
(  0.09442696,   0.04895706)
( 13.85808357,   0.03395526)
(  3.04491507,   0.03179233)
( 31.30293347,   0.03185319)
(  0.05128732,   0.02376444)
(  0.09635419,   0.02169466)
(  0.17853730,   0.01750213)
(  0.12152558,   0.03015192)
(  0.03221284,   0.03216546)
(  0.10986372,   0.02939943)
( 100.57686434,   0.04998021)
(  4.54662320,   0.06264639)
( 296.18124019,   0.04143467)
(  0.09020749,   0.02422864)
(  0.16052942,   0.00809136)
(  0.26785947,   0.02275230)
(  0.07616874,   0.03839638)
(  0.03825135,   0.01248658)
(  0.11000694,   0.05494363)};\label{line:deim_v_romnn_train}

\addplot [red, only marks, mark=o, mark size=1]
coordinates {
( 100.39017105,   0.01878609)
(  3.49995006,   0.01053297)
(  0.12504073,   0.03657047)
(  0.16536262,   0.03722634)
(  0.21202161,   0.02248219)
(  0.26609597,   0.02564125)
(  0.33098617,   0.04125064)
(  0.02848414,   0.02731054)
(  0.06175166,   0.01547970)
(  0.06607883,   0.01182947)
(  0.02917107,   0.02563175)
(  0.02339086,   0.02083726)
(  0.03625808,   0.01836958)
(  0.04604062,   0.03722730)
(  0.10240324,   0.02629014)
(  0.04771954,   0.03437760)
(  0.37613675,   0.02395702)
( 26.47324255,   0.02335248)
( 272.14561447,   0.02906314)
(  3.19142981,   0.03021269)
(  0.16284692,   0.03615993)
(  0.21010345,   0.03220961)
(  0.26597118,   0.01500449)
(  0.33653393,   0.02863762)
(  0.44246398,   0.02639380)
(  0.03713751,   0.02677102)
(  0.04668786,   0.03157620)
(  0.06686030,   0.01358292)
(  0.09607048,   0.01608449)
(  0.13425070,   0.01585765)
(  0.05046435,   0.02165116)
(  0.02242988,   0.02882965)
(  0.01970694,   0.01494911)
(  0.02692676,   0.01161077)
(  0.02963014,   0.01852827)
(  0.15197523,   0.02474738)
(  0.07353253,   0.03491148)
(  0.02938017,   0.02899157)
(  0.05978416,   0.02624471)
(  0.10405421,   0.03083135)
( 625.82164134,   0.03123512)
(  2.13943889,   0.02295305)
(  0.20235874,   0.02808195)
(  0.25762336,   0.02600410)
(  0.32760898,   0.01378897)
(  0.43527504,   0.01708425)
(  0.74585294,   0.02470142)
(  0.06899434,   0.03318531)
(  0.13293890,   0.01243260)
(  0.04015923,   0.02499198)
(  0.01779544,   0.03144882)
(  0.01277419,   0.01722334)
(  0.01423584,   0.00757926)
(  0.01259970,   0.02784174)
(  0.05238421,   0.04142949)
(  0.06701172,   0.02084144)
( 22.39835457,   0.04247768)
(  0.71221632,   0.04586525)
( 538.72200824,   0.04346990)
(  5.79013402,   0.02824612)
( 23.27191509,   0.03509380)
(  0.24379595,   0.02257584)
(  0.31014821,   0.01521108)
(  0.40703103,   0.01475263)
(  0.65046772,   0.01573797)
(  3.29237512,   0.02652944)
(  0.06947303,   0.02418241)
(  0.09395785,   0.02269628)
(  0.12788353,   0.01249368)
(  0.17105421,   0.01237791)
(  0.22324582,   0.02278368)
(  0.03392836,   0.02723191)
(  0.01566037,   0.02320197)
(  0.00939005,   0.00959396)
(  0.01371368,   0.00874888)
(  0.02861829,   0.03948210)
(  0.09602114,   0.03338410)
(  0.03861615,   0.03380542)
(  0.03630223,   0.02456147)
(  0.07013752,   0.01919761)
(  0.11215412,   0.04065278)
(  6.51134553,   0.05986598)
( 175.86330031,   0.03684954)
(  0.28827476,   0.01811210)
(  0.37263943,   0.01572675)
(  0.53765643,   0.03042262)
(  2.11725704,   0.01849833)
(  7.28143051,   0.02478694)
(  0.12069401,   0.01661867)
(  0.20950812,   0.00919676)
(  0.03146367,   0.02859294)
(  0.01936222,   0.01795675)
(  0.02078584,   0.00887351)
(  0.03441856,   0.01044382)
(  0.06079244,   0.04383227)
(  0.03048717,   0.02936239)
(  0.06960167,   0.01777446)};\label{line:deim_v_romnn_test}

\end{groupplot}\end{tikzpicture}
 \caption{Comparison of the error in the (D)EIM and ROM-NN
          ($k_\ubm = 8$) with respect
          to the HDM solution when trained with $\Xibold_0^a$
          (\textit{left}) or $\Xibold_0^b$ (\textit{right}) for
          each point in $\Xibold^*$. The individual marks correspond to
          the (D)EIM error vs. the ROM error for training
          (\ref{line:deim_v_romnn_train}) and testing
          (\ref{line:deim_v_romnn_test}) points. All entries that
          lie below the line of identity (\ref{line:eye}) correspond
          to parameters where the ROM-NN is more accurate than (D)EIM.
          For both training cases, far more points lie below
          the line of identity indicating the ROM-NN is more accurate
          across the testing set $\Xibold^*$ than (D)EIM.}
 \label{fig:vburg_train_all_err0}
\end{figure}
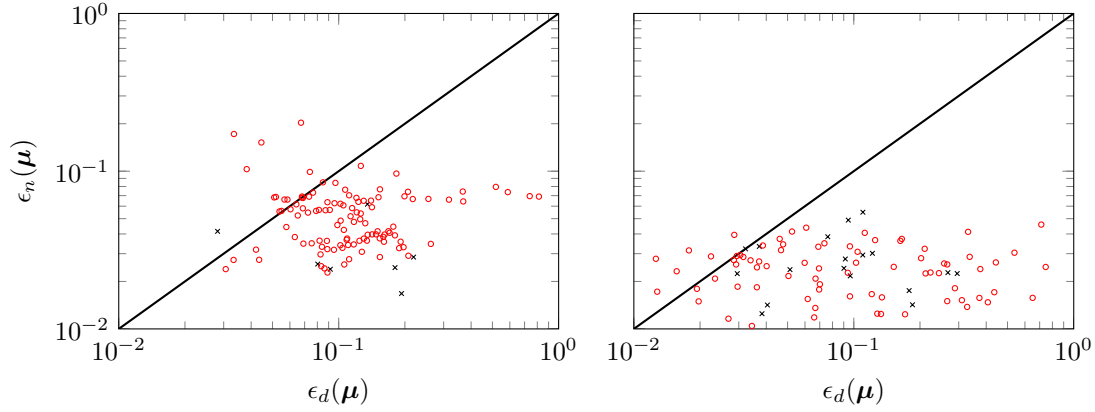
}

%\begin{itemize}
% \item For 2x2x2 training
% \begin{itemize}
%  \item Minimum error: comparable between (D)EIM and ROM-NN, both train/test
%  \item Maximum error: huge for (D)EIM b/c of unstable, about 1 order of
%        magintude larger for ROM-NN than ROM, but $6\%$ for training and
%        $22\%$ for testing, unacceptably large
%  \item Median error: factor of 5 larger than ROM and factor of 5 smaller
%        than (D)EIM; less than $3\%$ for training and about $5\%$ for
%        testing
% \end{itemize}
% \item For 3x3x3 training
% \begin{itemize}
%  \item (D)EIM less stable: more unstable points, both training/testing;
%        stable points a bit more accurate (median errors decrease)
%  \item ROM-NN: training errors similar, testing errors smaller;
%        additional training leads to better prediction
%  \item ROM-NN: median errors less than $3\%$ (both), max errors
%        about $6\%$
% \end{itemize}
%\end{itemize}

\subsection{Two-dimensional premixed H$_2$-air flame model}
\label{sec:numexp:2d}
The second numerical experiment we consider is solution of a
simplified model of a premixed H$_2$-air flame at a constant and
uniform pressure, in a constant, divergence-free velocity field, and
with constant, uniform diffusivities for all species and
temperature in the domain $\Omega \coloneqq [0, L_x] \times [0, L_y]$,
where $L_x = 18$mm and $L_y = 9$mm, over the time interval
$[0, T]$, $T = 0.06$s. The one-step reaction mechanism is
$2\text{H}_2 + \text{O}_2 \rightarrow 2\text{H}_2\text{O}$.
The PDE model of this system \cite{buffoni2010projection} is
\begin{equation}
 \begin{aligned}
  \partial_t U(x, t, \mubold) -
  \kappa \Delta U(x, t, \mubold) +
  \beta \cdot \nabla U(x, t, \mubold) &=
  \Ncal(U, \mubold),
  &&x \in \Omega, &&t \in [0, T], &&\mubold \in \Dcal, \\
  U(x, t, \mubold) &= U_D(x),
  &&x \in \Gamma_D,
  && t \in [0, T], && \mubold \in \Dcal, \\
  \nabla U(x, t, \mubold) \cdot n(x) &= 0,
  &&x \in \Gamma_N,
  && t \in [0, T], && \mubold \in \Dcal, \\
  U(x, 0, \mubold) &= U_0, && x \in \Omega, && && \mubold \in \Dcal
 \end{aligned}
\end{equation}
with solution
\begin{equation}
 \func{U}{\Omega \times [0, T]\times\Dcal}{\Rbb^4}, \qquad
 (x, t, \mubold) \mapsto
 \begin{bmatrix}
  Y_F(x, t, \mubold) \\
  Y_O(x, t, \mubold) \\
  Y_P(x, t, \mubold) \\
  \Theta(x, t, \mubold)
 \end{bmatrix},
\end{equation}
where $\func{n}{\partial\Omega}{\Rbb^2}$ is the outward unit normal,
$\func{Y_i}{\Omega \times [0, T] \times \Dcal}{\Rbb}$ is the
mass fraction of the hydrogen fuel ($i = F$), oxygen ($i = O$), and
water product ($i = P$), and
$\func{\Theta}{\Omega \times [0, T] \times \Dcal}{\Rbb}$
is the temperature. The domain boundary is split into six
segments (Figure~\ref{fig:hcube0})
\begin{equation}
 \partial\Omega = \bigcup_{i=1}^6 \Gamma_i, \qquad
 \Gamma_D \coloneqq \bigcup_{i=1}^3 \Gamma_i, \qquad
 \Gamma_N \coloneqq \bigcup_{i=4}^6 \Gamma_i
\end{equation}
with the following essential boundary conditions prescribed on
$\Gamma_D\subset\partial\Omega$
\begin{equation}
 \func{U_D}{\Omega}{\Rbb^4}, \qquad
 x \mapsto
 \begin{cases}
  (0, 0, 0, 300) & x \in \Gamma_1\cup\Gamma_3 \\
  (0.0282, 0.2259, 0, 950) & x \in \Gamma_2
 \end{cases}
\end{equation}
and homogeneous natural boundary conditions on $\Gamma_N\subset\partial\Omega$.
\ifbool{fastcompile}{}{
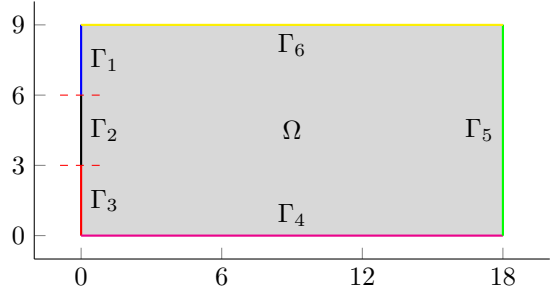
\begin{figure}
 \centering
 \begin{tikzpicture}
\begin{axis}[
axis equal image,
axis x line*=bottom,
axis y line*=left,
xtick={0, 6, 12, 18},
ytick={0, 3, 6, 9},
ymax=10,
xmax=20,
xmin=-2,
ymin=-1]
\addplot [black, solid, opacity=0.6, fill=lightgray, forget plot]
coordinates {
(  0.00000000,   0.00000000)
( 18.00000000,   0.00000000)
( 18.00000000,   9.00000000)
(  0.00000000,   9.00000000)
(  0.00000000,   0.00000000)};

\addplot [red, thick, forget plot]
coordinates {
(  0.00000000,   0.00000000)
(  0.00000000,   3.00000000)};

\addplot [black, thick, forget plot]
coordinates {
(  0.00000000,   3.00000000)
(  0.00000000,   6.00000000)};

\addplot [blue, thick, forget plot]
coordinates {
(  0.00000000,   6.00000000)
(  0.00000000,   9.00000000)};

\addplot [magenta, thick, forget plot]
coordinates {
(  0.00000000,   0.00000000)
( 18.00000000,   0.00000000)};

\addplot [green, thick, forget plot]
coordinates {
( 18.00000000,   0.00000000)
( 18.00000000,   9.00000000)};

\addplot [yellow, thick, forget plot]
coordinates {
(  0.00000000,   9.00000000)
( 18.00000000,   9.00000000)};

\addplot [red, dashed, forget plot]
coordinates {
( -0.90000000,   3.00000000)
(  0.90000000,   3.00000000)};

\addplot [red, dashed, forget plot]
coordinates {
( -0.90000000,   6.00000000)
(  0.90000000,   6.00000000)};

\node[color=black]    at    (axis cs:9.0, 4.5) {$\Omega$};
\node[right, color=black]    at    (axis cs:0, 7.5) {$\Gamma_1$};
\node[right, color=black]    at    (axis cs:0, 4.5) {$\Gamma_2$};
\node[right, color=black]    at    (axis cs:0, 1.5) {$\Gamma_3$};
\node[above, color=black]    at    (axis cs:9.0, 0) {$\Gamma_4$};
\node[left, color=black]    at    (axis cs:18, 4.5) {$\Gamma_5$};
\node[below, color=black]    at    (axis cs:9.0, 9) {$\Gamma_6$};
\end{axis}
\end{tikzpicture}
 \caption{Schematic setup for the hydrogen-air flame (units: mm).}
 \label{fig:hcube0}
\end{figure}
}
The nonlinear reaction source term is of Arrhenius type and modeled as in
Cuenot and Poinsot \cite{cuenot1996asymptotic} as
$\Ncal(U, \mubold) =
  [\Ncal_F(U, \mubold),
   \Ncal_O(U, \mubold),
   \Ncal_P(U, \mubold),
   \Ncal_\Theta(U, \mubold)]$, where
\begin{equation}
 \begin{aligned}
  \Ncal_i(U, \mubold) &= -\nu_i \left(\frac{W_i}{\rho}\right)
                          \left(\frac{\rho Y_F}{W_F}\right)^{\nu_F}
                          \left(\frac{\rho Y_O}{W_O}\right)^{\nu_O}
                          A(\mubold)
                     \exp\left(-\frac{E(\mubold)}{R\Theta}\right) \\
  \Ncal_\Theta(U, \mubold) &= \Ncal_P(U, \mubold) Q
 \end{aligned}
\end{equation}
for $i = F, O, P$.
The divergence-free velocity field is $\beta = (50, 0)$ cm/sec.
The diffusivities are $\kappa = 2.0$ cm$^2$/sec and the density of
the mixture is $\rho = 1.39\times 10^{-3}$ gr/cm$^3$. The molecular
weights are $W_F = 2.016$, $W_O = 31.9$, $W_P = 18$ gr/mol,
the stoichiometric coefficients are $\nu_F = 2$, $\nu_O = 1$,
$\nu_P = 2$, the heat of reaction is $Q = 9800$K, and the universal
gas constant is $R = 8.314472$ J/(mol K). The parameter space is
taken as
$\Dcal \coloneqq [2.3375\times 10^{12}, 6.2\times 10^{12}]
                 \times [5625.5, 9000]$.
For any
$\mubold \in \Dcal$ where $\mubold = (\mu_1, \mu_2)$, the
parametrized pre-exponential factor and activation energy are
taken as
\begin{equation}
 A(\mubold) = \mu_1, \quad E(\mubold) = \mu_2.
\end{equation}
At $t = 0$, the domain is considered empty at a temperature of $300K$, i.e.,
$U_0 = (0, 0, 0, 300)$. The PDE is discretized in space using the finite
difference method on a grid of $40 \times 20$ with essential boundary
conditions strongly enforced to yield a dynamical system of the form
(\ref{eqn:semidisc}) with a total of $N_\ubm = 2736$ spatial degrees of
freedom. The dynamical system is discretized in time using a two-stage
diagonally implicit Runge-Kutta method with $50$ time steps; see
Figure~\ref{fig:adr_stvc} for solution snapshots for a representative
parameter configuration.

\ifbool{fastcompile}{}{
\begin{figure}
 \centering
 \begin{tikzpicture}
\begin{groupplot}[
  group style={
      group size=2 by 3,
      horizontal sep=5.0mm,
      vertical sep=5.0mm,
  },
  width=0.5\textwidth,
  axis equal image,
  axis line style={gray},
  axis x line*=bottom,
  axis y line*=left, 
  xtick = {0, 6, 12, 18},
  ytick = {0, 3, 6, 9},
  xmin=0.0, xmax=18,
  ymin=0.0, ymax=9,
]

\nextgroupplot[ylabel={$x_2$ (mm)}, xticklabels={,,}]
\addplot graphics [xmin=0.0, xmax=18, ymin=0.0, ymax=9] {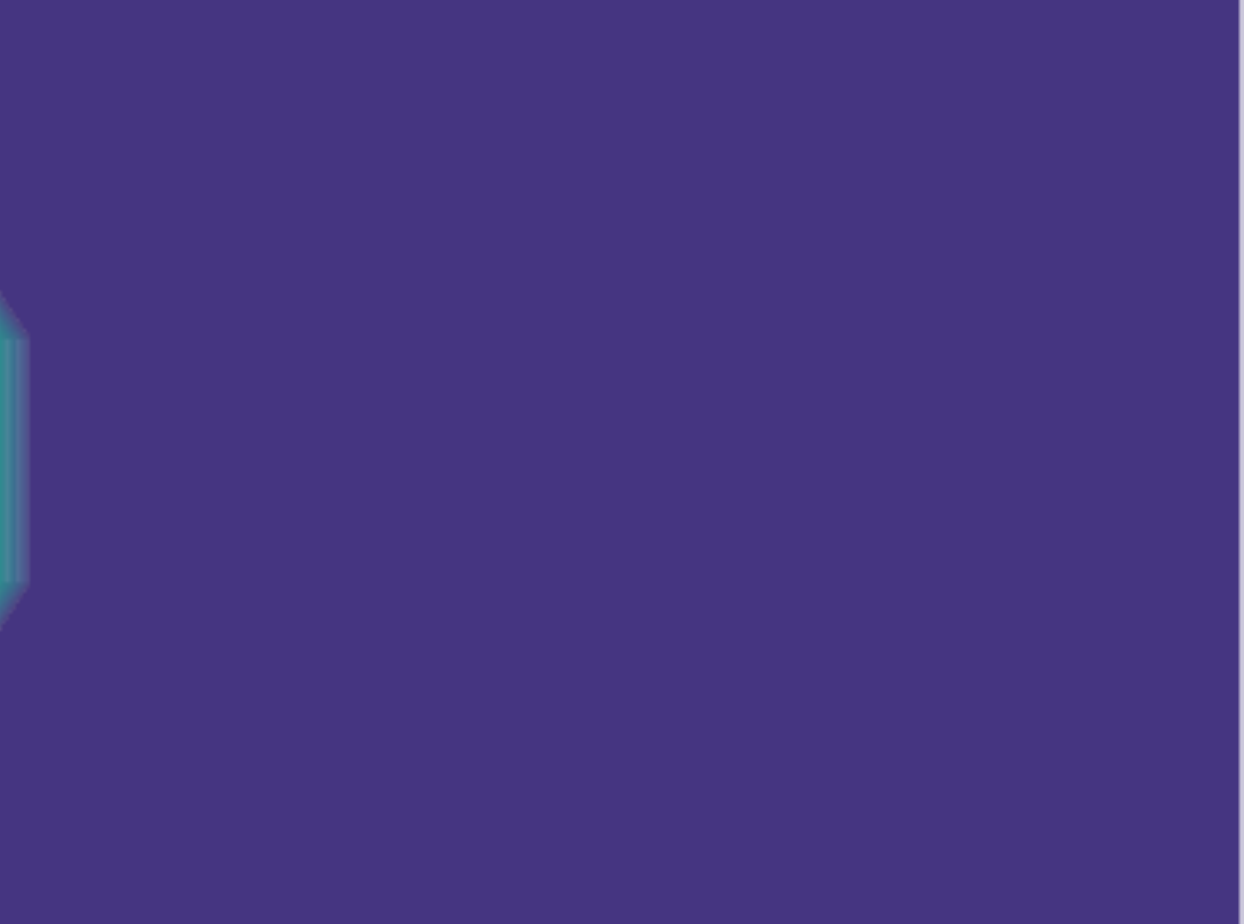};

\nextgroupplot[xticklabels={,,}, yticklabels={,,}]
\addplot graphics [xmin=0.0, xmax=18, ymin=0.0, ymax=9] {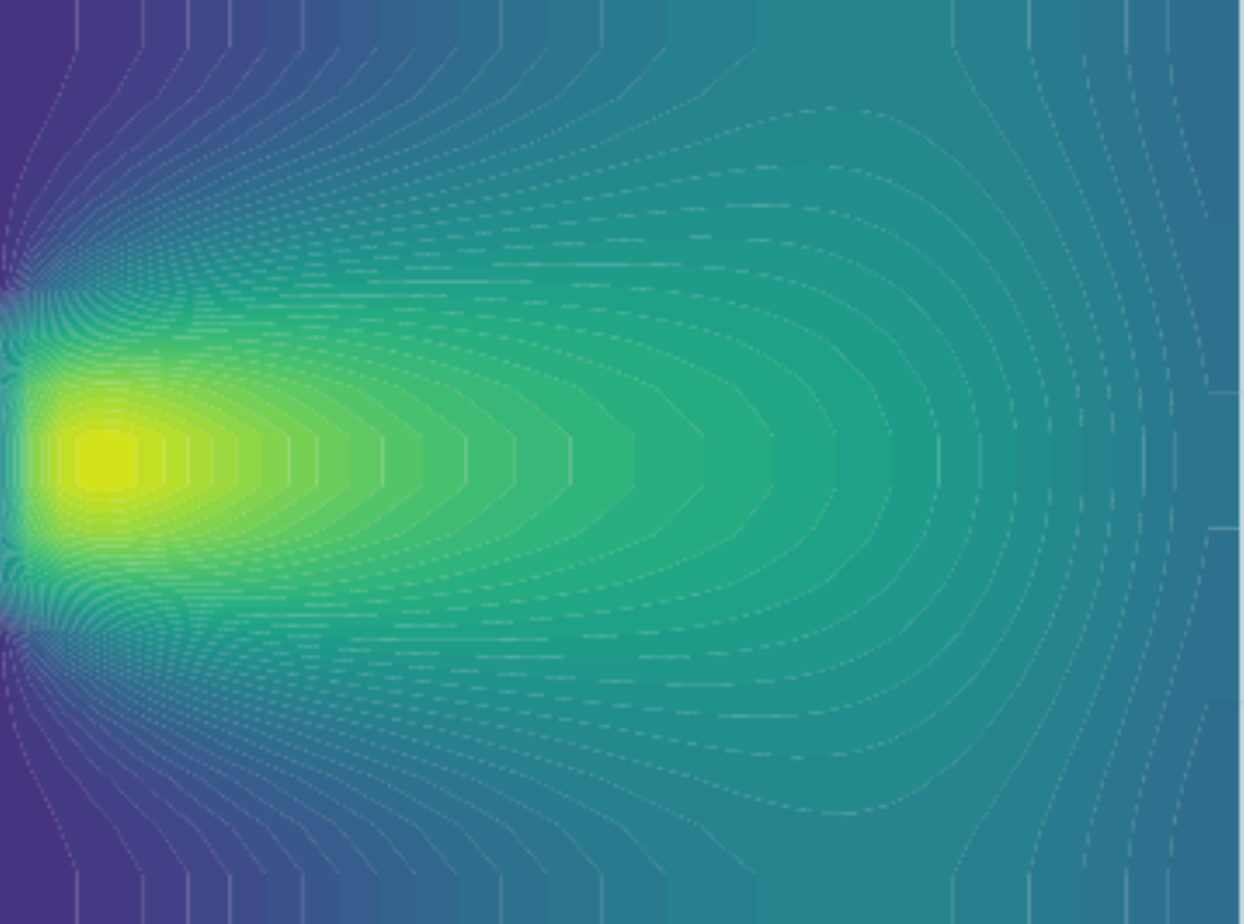};

\nextgroupplot[ylabel={$x_2$ (mm)}, xticklabels={,,}]
\addplot graphics [xmin=0.0, xmax=18, ymin=0.0, ymax=9] {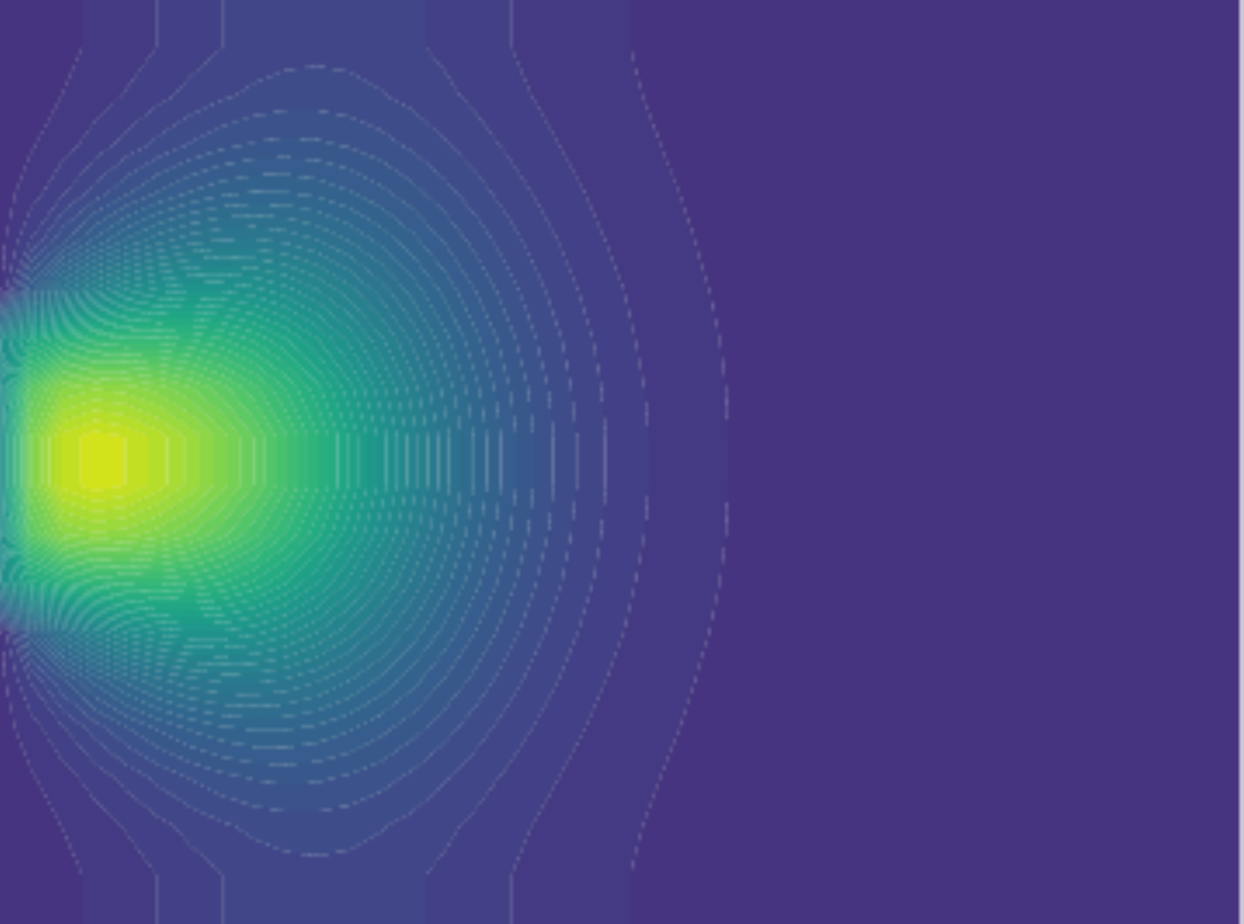};

\nextgroupplot[xticklabels={,,}, yticklabels={,,}]
\addplot graphics [xmin=0.0, xmax=18, ymin=0.0, ymax=9] {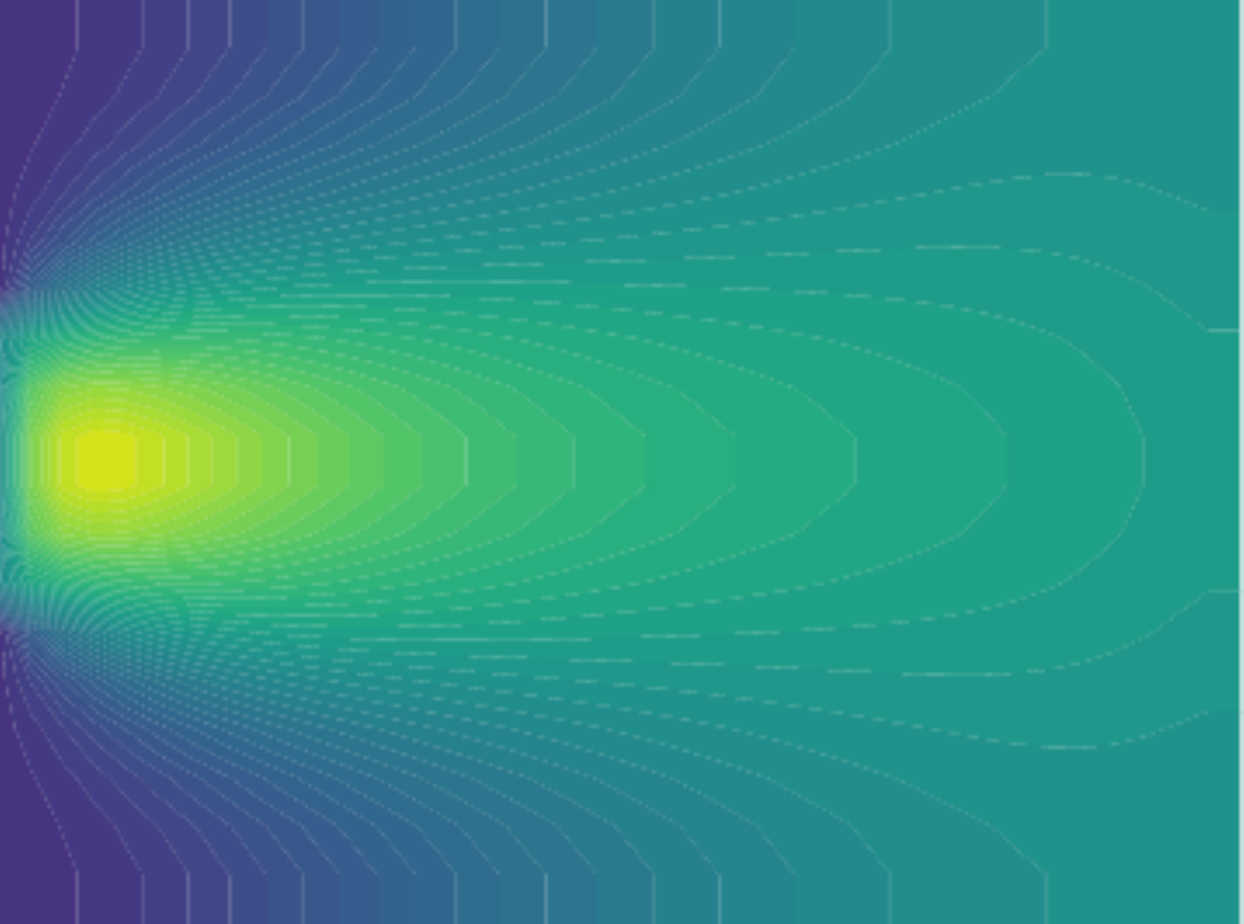};

\nextgroupplot[xlabel={$x_1$ (mm)}, ylabel={$x_2$ (mm)}]
\addplot graphics [xmin=0.0, xmax=18, ymin=0.0, ymax=9] {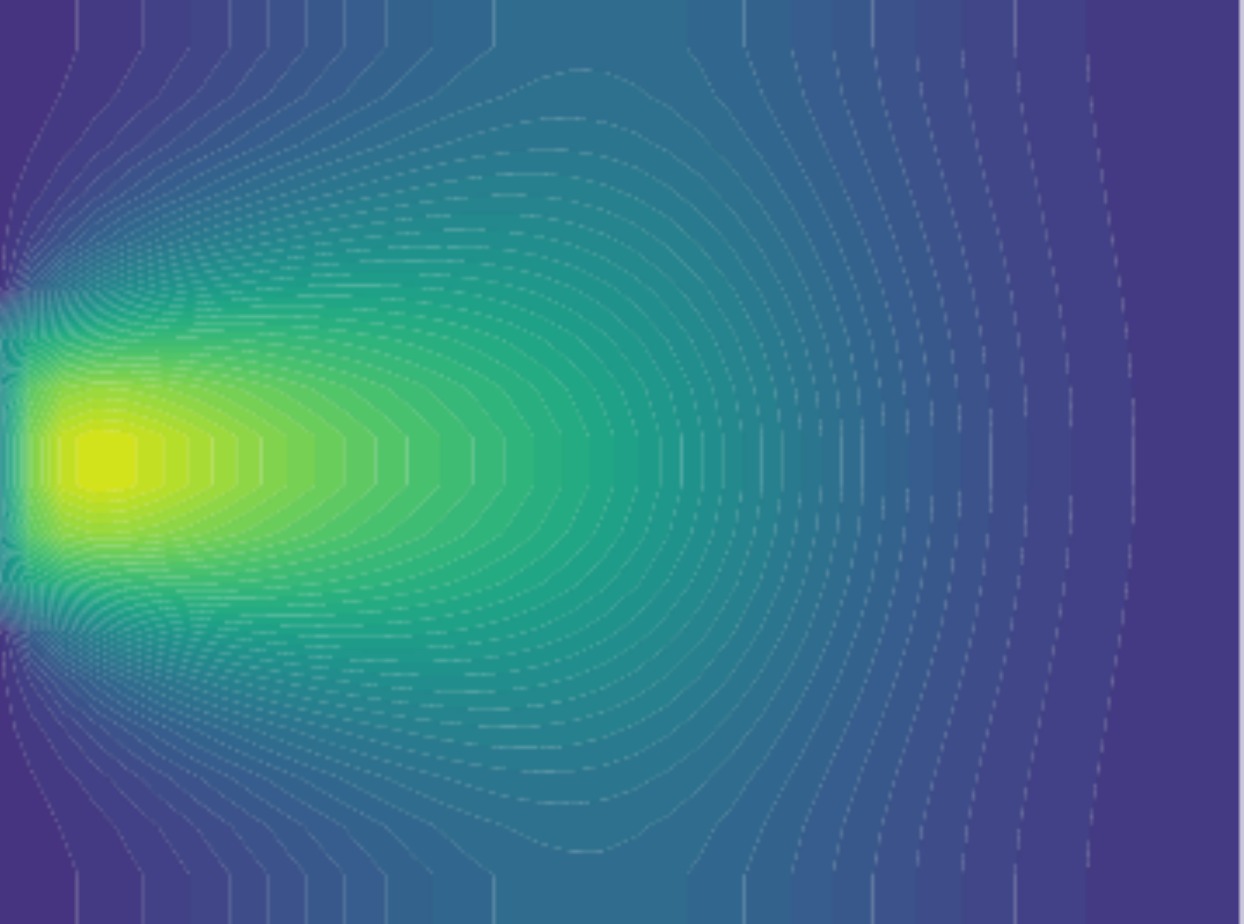};

\nextgroupplot[xlabel={$x_1$ (mm)}, yticklabels={,,}]
\addplot graphics [xmin=0.0, xmax=18, ymin=0.0, ymax=9] {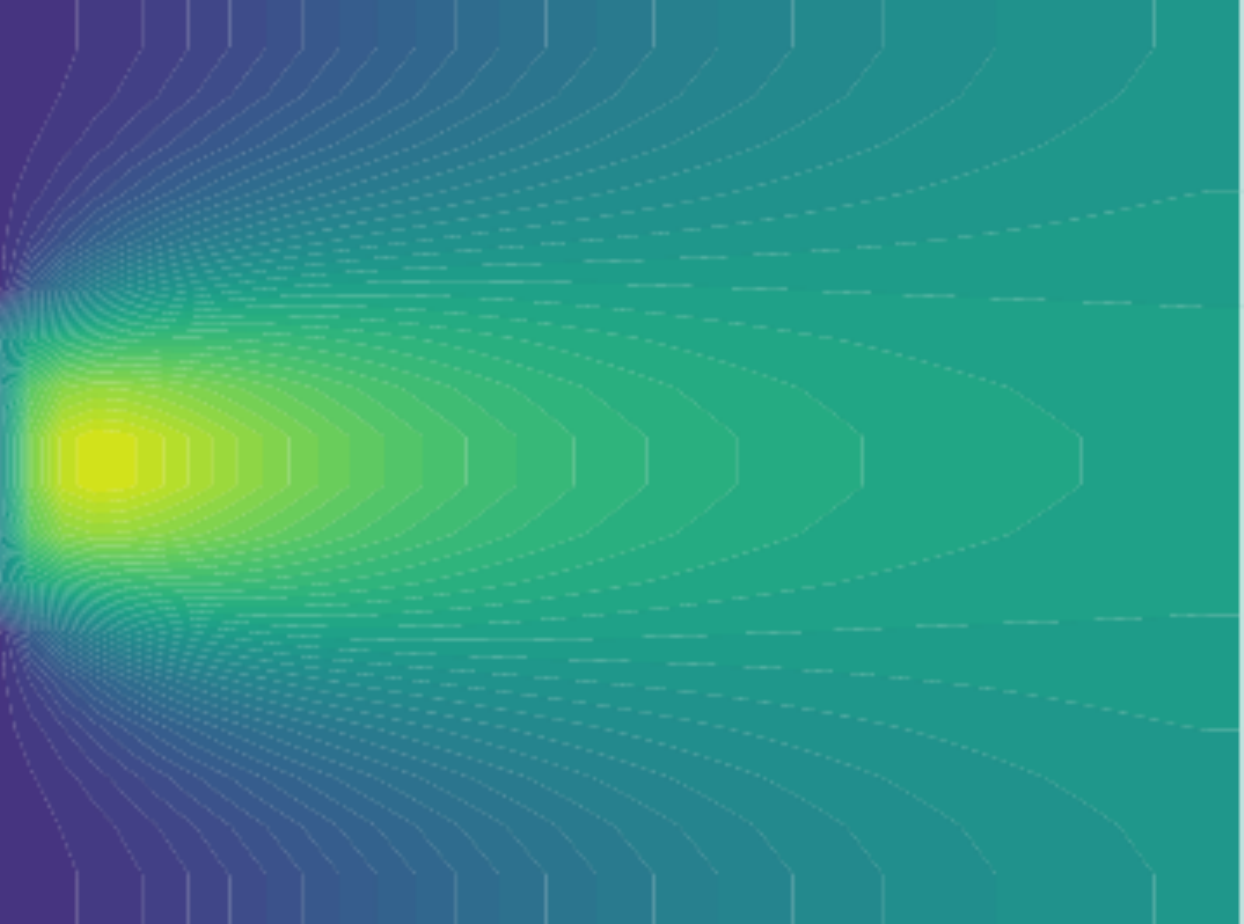};
\end{groupplot}

\end{tikzpicture}
 \caption{Snapshots of advection-diffusion-reaction temperature
          ($\Theta(x,t,\mubold)$) field
          ($t = 0, 0.012, 0.024, 0.036, 0.048, 0.06$; top-to-bottom then
          left-to-right) at $\mubold = (5.1125\times 10^{12}, 6187.91667)$
          computed using the HDM.}
 \label{fig:adr_stvc}
\end{figure}
}

For this problem, we define the testing set $\Xibold^*$ as the uniform
sampling of $\Dcal$ on a $7 \times 7$ grid for a total of
$|\Xibold^*| = 49$ test points. The training set is a uniform
sampling of $\Dcal$ on a $4\times 4$ grid, respectively.
By construction, $\Xibold_0 \subset \Xibold^*$.
Similar to the previous section, we construct a POD-Galerkin
ROM without hyperreduction, one accelerated with (D)EIM, and
one accelerated with the neural network for a range of basis
sizes ($k_\ubm$) and subsequently test each model on all points
in $\Xibold^*$.

The reduced-order model without hyperreduction is the most stable and
accurate method and demonstrates deep convergence under refinement of
$k_\ubm$, even when error metric is aggregated over both testing and
training points. For this problem, (D)EIM was highly unstable; a
basis of size $k_\ubm = 150$ was required for stability, while the
other methods were stable for a basis of size $k_\ubm = 80$. When
(D)EIM is stable, it is more accurate than the proposed ROM-NN method.
The ROM-NN method is stable for all basis sizes considered, but does not
exhibit deep convergence
(the error is about $2\%$ for basis sizes $k_\ubm \in [80, 200]$).
We have run extensive numerical tests to confirm this is due to the scaling
of the entries of $\fbm_r$, which can vary by up to seven orders of magnitude
for this problem. This implies the loss function used to train the FCNN is
heavily biased toward the important coefficients and the smallest coefficients
are all but ignored in the training. As a result, the relative error of the
small coefficients (required for deep convergence) is large. Weighting the
loss function to offset the massively different scale of the loss
function terms causes the approximation of small coefficients to improve, but
the accuracy of the large (important) coefficients degrades for the networks
considered in this work, which causes the overall error to increase.
Because the ROM-NN does not exhibit deep convergence, there is little
point in refining the basis beyond $k_\ubm = 80$ because the reduction
in error is negligible. The maximum and median errors of the ROM-NN
solution ($k_\ubm = 80$) across all training and testing parameters are
small ($2-3\%$) (Figure~\ref{fig:adr_train0_err_vs_ny}).

%, which is expected since $\hat\fbm$ has more
%inputs/outputs as $k_\ubm$ increases and becomes more difficult to train.

\ifbool{fastcompile}{}{
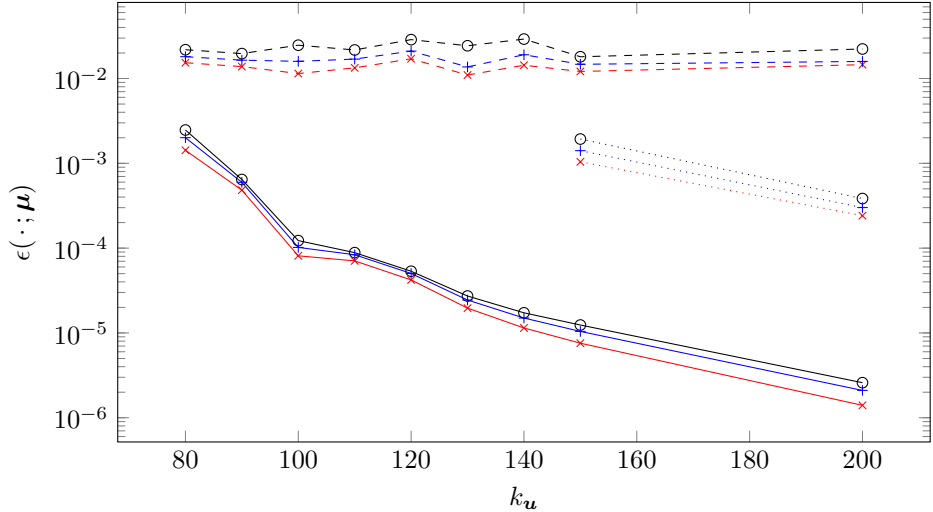
\begin{figure}
 \centering
 \begin{tikzpicture}
\begin{axis}[
height=0.45\textwidth,
width=0.75\textwidth,
ymode=log,
xlabel={$k_\ubm$},
ylabel={$\epsilon(\,\cdot\,;\mubold)$}]
\addplot [solid, black, mark options={solid}, mark=o, mark size=2]
coordinates {
( 80.00000000,   0.00247058)
( 90.00000000,   0.00064477)
( 100.00000000,   0.00012264)
( 110.00000000,   0.00008847)
( 120.00000000,   0.00005338)
( 130.00000000,   0.00002724)
( 140.00000000,   0.00001727)
( 150.00000000,   0.00001240)
( 200.00000000,   0.00000259)};\label{line:adr:rom_max}

\addplot [solid, red, mark options={solid}, mark=x, mark size=2]
coordinates {
( 80.00000000,   0.00142507)
( 90.00000000,   0.00048240)
( 100.00000000,   0.00008116)
( 110.00000000,   0.00007072)
( 120.00000000,   0.00004219)
( 130.00000000,   0.00001961)
( 140.00000000,   0.00001144)
( 150.00000000,   0.00000757)
( 200.00000000,   0.00000140)};\label{line:adr:rom_min}

\addplot [solid, blue, mark options={solid}, mark=+, mark size=2]
coordinates {
( 80.00000000,   0.00200580)
( 90.00000000,   0.00059754)
( 100.00000000,   0.00010198)
( 110.00000000,   0.00008356)
( 120.00000000,   0.00005025)
( 130.00000000,   0.00002443)
( 140.00000000,   0.00001499)
( 150.00000000,   0.00001042)
( 200.00000000,   0.00000210)};\label{line:adr:rom_med}

\addplot [dashed, black, mark options={solid}, mark=o, mark size=2]
coordinates {
( 80.00000000,   0.02185277)
( 90.00000000,   0.01967411)
( 100.00000000,   0.02470669)
( 110.00000000,   0.02170276)
( 120.00000000,   0.02874913)
( 130.00000000,   0.02433513)
( 140.00000000,   0.02918716)
( 150.00000000,   0.01805090)
( 200.00000000,   0.02229232)};\label{line:adr:romnn_max}

\addplot [dashed, red, mark options={solid}, mark=x, mark size=2]
coordinates {
( 80.00000000,   0.01535441)
( 90.00000000,   0.01381741)
( 100.00000000,   0.01144429)
( 110.00000000,   0.01335115)
( 120.00000000,   0.01708672)
( 130.00000000,   0.01099225)
( 140.00000000,   0.01433767)
( 150.00000000,   0.01211177)
( 200.00000000,   0.01456074)};\label{line:adr:romnn_min}

\addplot [dashed, blue, mark options={solid}, mark=+, mark size=2]
coordinates {
( 80.00000000,   0.01809733)
( 90.00000000,   0.01645671)
( 100.00000000,   0.01594476)
( 110.00000000,   0.01694089)
( 120.00000000,   0.02104021)
( 130.00000000,   0.01369702)
( 140.00000000,   0.01907132)
( 150.00000000,   0.01471593)
( 200.00000000,   0.01590362)};\label{line:adr:romnn_med}

\addplot [dotted, black, mark options={solid}, mark=o, mark size=2]
coordinates {
( 80.00000000,          nan)
( 90.00000000,          nan)
( 100.00000000,          nan)
( 110.00000000,          nan)
( 120.00000000,          nan)
( 130.00000000,          nan)
( 140.00000000,          nan)
( 150.00000000,   0.00193370)
( 200.00000000,   0.00038429)};\label{line:adr:deim_max}

\addplot [dotted, red, mark options={solid}, mark=x, mark size=2]
coordinates {
( 80.00000000,          nan)
( 90.00000000,          nan)
( 100.00000000,          nan)
( 110.00000000,          nan)
( 120.00000000,          nan)
( 130.00000000,          nan)
( 140.00000000,          nan)
( 150.00000000,   0.00104154)
( 200.00000000,   0.00024143)};\label{line:adr:deim_min}

\addplot [dotted, blue, mark options={solid}, mark=+, mark size=2]
coordinates {
( 80.00000000,          nan)
( 90.00000000,          nan)
( 100.00000000,          nan)
( 110.00000000,          nan)
( 120.00000000,          nan)
( 130.00000000,          nan)
( 140.00000000,          nan)
( 150.00000000,   0.00140979)
( 200.00000000,   0.00030079)};\label{line:adr:deim_med}

\end{axis}
\end{tikzpicture}
 \caption{The maximum (circle), minimum (cross), and median (plus)
          error of the ROM (solid), (D)EIM (dotted), and ROM-NN (dashed)
          over the set $\Xibold^*$. While the ROM error statistics demonstrate
          deep convergence under refinement of $k_\ubm$, the ROM-NN does not
          (expected because $\hat\fbm$ has more inputs/outputs as $k_\ubm$
          increases and therefore becomes more difficult to train). Even
          though the ROM-NN does not have deep convergence, its maximum
          and median errors are small ($3\%$). (D)EIM
          is unstable for $k_\ubm < 150$; however, once the basis is
          sufficiently large $k_\ubm \geq 150$ for stability, it is more
          accurate than the ROM-NN.
          Legend: maximum ROM error (\ref{line:adr:rom_max}),
                  minimum ROM error (\ref{line:adr:rom_min}),
                  median  ROM error (\ref{line:adr:rom_med}),
                  maximum (D)EIM error (\ref{line:adr:deim_max}),
                  minimum (D)EIM error (\ref{line:adr:deim_min}),
                  median  (D)EIM error (\ref{line:adr:deim_med}),
                  maximum ROM-NN error (\ref{line:adr:romnn_max}),
                  minimum ROM-NN error (\ref{line:adr:romnn_min}),
                  median  ROM-NN error (\ref{line:adr:romnn_med}).}
 \label{fig:adr_train0_err_vs_ny}
\end{figure}
}

\section{Conclusion}
\label{sec:conclude}
This work developed a \textit{non-intrusive} acceleration technique
for projection-based model reduction of nonlinear dynamical systems using
deep neural networks. The approach is non-intrusive in the sense that
it treats both the original dynamical system and neural network code as
black boxes. The method is trained using the same HDM solutions computed
to train the reduced basis, i.e., no new simulations are required to
train the neural network approximation of the ROM velocity function,
only evaluations of the ROM velocity at existing snapshots. Unlike
traditional hyperreduction methods, this does not require modification
of the underlying dynamical system code because, once the neural
network is trained, only forward propagation (and symbolic differentiation)
through the network is required to compute the approximate velocity function
and its derivative. Aside from the benefit of non-intrusivity, the proposed
method is more stable and accurate at both training and testing points
in the limit of a small basis than the popular (D)EIM hyperreduction
for the two dynamical systems considered (semi-discretizations of
nonlinear, hyperbolic PDEs). Given we used uniform sampling to train
the ROM-NN, this approach may be most appropriate for problems with a
limited number of parameters. Future work will explore whether the
amount of training can be reduced using POD-Greedy sampling
\cite{haasdonk2013convergence} without sacrificing stability or
parametric robustness. We will also perform a careful study of the
computational cost of the proposed approach and the benefits of using
GPUs to train and pass through the neural network in future work.

\bibliographystyle{elsarticle-num}
\bibliography{biblio}

\end{document}